\documentclass{amsart}
\usepackage{amssymb}
\usepackage{amsbsy}
\usepackage{amscd}
\usepackage{epic}
\usepackage{eepic}
\usepackage[dvips]{graphics}
\DeclareGraphicsExtensions{.eps}
\makeatletter
\def\Bbb{\mathbb}
\def\frak{\mathfrak}

\newenvironment{pf*}[1]{\proof[#1]}{\endproof}
\newcommand{\rom}{\textup}
\makeatother
\newenvironment{aenume}{%
  \begin{enumerate}%
  }{\end{enumerate}}
\newtheorem{Theorem}[equation]{Theorem}
\newtheorem{Corollary}[equation]{Corollary}
\newtheorem{Lemma}[equation]{Lemma}
\newtheorem{Proposition}[equation]{Proposition}

\theoremstyle{definition}
\newtheorem{Definition}[equation]{Definition}

\theoremstyle{remark}
\newtheorem{Remark}[equation]{Remark}

\numberwithin{equation}{section}
\numberwithin{figure}{section}

\newcommand{\thmref}[1]{Theorem~\ref{#1}}
\newcommand{\secref}[1]{\S\ref{#1}}
\newcommand{\lemref}[1]{Lemma~\ref{#1}}
\newcommand{\propref}[1]{Proposition~\ref{#1}}
\newcommand{\corref}[1]{Corollary~\ref{#1}}

\newcommand{\defeq}{\overset{\operatorname{\scriptstyle def.}}{=}}
\newcommand{\C}{{\Bbb C}}

\newcommand{\Q}{{\Bbb Q}}

\newcommand{\proj}{{\Bbb P}}
\newcommand{\CP}{\operatorname{\C P}}

\newcommand{\SL}{\operatorname{\rm SL}}

\newcommand{\GL}{\operatorname{GL}}

\newcommand{\Hom}{\operatorname{Hom}}

\newcommand{\Ker}{\operatorname{Ker}}

\newcommand{\Ima}{\operatorname{Im}}

\newcommand{\rank}{\operatorname{rank}}

\newcommand{\tr}{\operatorname{tr}}

\newcommand{\tautR}{\mathcal R} 
\newcommand{\Hilb}[2]{\operatorname{Hilb}^{#2}(#1)} 

\newcommand{\shfO}{\mathcal O} 
\newcommand{\uZ}{\mathcal Z} 
\newcommand{\Supp}{\operatorname{Supp}} 
\newcommand{\ch}{\operatorname{ch}} 
\begin{document}
\title{McKay correspondence and Hilbert schemes
in dimension three}
\author{Yukari Ito}
\address{Department of Mathematics, Tokyo Metropolitan University, Hachioji, 
Tokyo 192-03, Japan}
\curraddr{Lehrestuhl f\"ur Mathematik VI, Universit\"at Mannheim, 68131 
Mannheim, Germany}
\email{yukari@math.metro-u.ac.jp}
\author{Hiraku Nakajima}
\address{Department of Mathematics, Kyoto University, Kyoto 606-01, Japan}
\email{nakajima@kusm.kyoto-u.ac.jp}
\thanks{The first author is partially supported by the Grant-in-aid for 
Scientific Research (No.09740036), the Ministry of Education and the F\^ujukai 
Foundation.
The second author is partially supported by the Grant-in-aid
for Scientific Research (No.09740008), the Ministry of Education, Japan.}
\begin{abstract}
Let $G$ be a nontrivial finite subgroup of $\SL_n(\C)$. Suppose that
the quotient singularity $\C^n/G$ has a crepant resolution $\pi\colon
X\to \C^n/G$ (i.e.\ $K_X = \shfO_X$). There is a slightly imprecise
conjecture, called the McKay correspondence, stating that there is a
relation between the Grothendieck group (or (co)homology group) of $X$
and the representations (or conjugacy classes) of $G$ with a ``certain
compatibility'' between the intersection product and the tensor
product (see e.g. \cite{Maizuru}). The purpose of this paper is to
give more precise formulation of the conjecture when $X$ can be given
as a certain variety associated with the Hilbert scheme of points in
$\C^n$. We give the proof of this new conjecture for an abelian
subgroup $G$ of $\SL_3(\C)$.
\end{abstract}

\maketitle

\section{Introduction}\label{sec:intro}

Let $G$ be a nontrivial finite subgroup of $\SL_n(\C)$ and let $X$ be
the scheme parametrising $0$-dimensional subschemes $Z$ of $\C^n$
satisfying the following three conditions:
\begin{enumerate}
\item the length of $Z$ is equal to $\# G = \text{the order of $G$}$.
\item $Z$ is invariant under the $G$-action.
\item $H^0(\shfO_Z)$ is the regular representation of $G$.
\end{enumerate}
This is a union of components (possibly one component) of fixed points
of the $G$-action on the Hilbert scheme of $\# G$-points in
$\C^n$. (See \cite{Lecture} for survey on Hilbert schemes of points.)
If $Z$ consists of pairwise distinct points, the above conditions mean
that $Z$ is a single $G$-orbit. We have a natural morphism (the
Hilbert-Chow morphism) $\pi\colon X\to \C^n/G$ which is a crepant
resolution of singularities under certain assumption on $G$ as
explained later.

In $n = 2$ (i.e.\  simple singularities), it was observed by
Ginzburg-Kapranov \cite{GK} and Ito-Nakamura \cite{IN} that $X$ is
nonsingular, and in fact, is the minimal resolution of $\C^2/G$ (see
also \cite[Chapter~4]{Lecture}). The proof is based on the fact that
Hilbert schemes of points on $\C^2$ are nonsingular symplectic
manifolds. (See e.g., \cite[Chapter~1]{Lecture}.) Since Hilbert
schemes are singular in higher dimensions in general, the proof is not
applied to the $3$-dimensional case. Hence we were surprised when
Nakamura proved that $X$ is nonsingular when $G$ is an abelian
subgroup of $\SL_3(\C)$ \cite{Iku}. He conjectured the same is true
for any $G\subset\SL_3(\C)$, and it is still an open problem. Anyhow,
it seems reasonable to consider $X$ as a first candidate for crepant
resolutions of $\C^n/G$.

Now we explain the McKay correspondence. We first recall the
$2$-dimensional situation. As $X$ is the minimal resolution of
$\C^2/G$, it is well-known that the exceptional set consists of
projective lines intersecting transversely. Let us denote by $C_k$ the
irreducible component. The intersection matrix $C_k\cdot C_l$ is given
by the negative of the Cartan matrix.
On the other hand, McKay~\cite{Mc} considered the irreducible
representations of $G$ and the decomposition of a tensor product 
\begin{equation*}
  Q \otimes \rho_l = \bigoplus_k a_{kl} \rho_k,
\end{equation*}
where $\{\rho_k\}_{k=0}^r$ be the set of isomorphism classes of
irreducible representations of $G$ and $Q$ is the $2$-dimensional
representation given by the inclusion $G\subset\SL_2(\C)$. Then he
observed that $(2 \delta_{kl} - a_{kl})$ is the extended Cartan
matrix. The trivial representation, denoted by $\rho_0$, corresponds
to the extra entry added to the finite Cartan matrix, which turns out
to be the same as realized by the intersection matrix.
The correspondences are summarized as follows:
\begin{center}
\begin{tabular}[t]{lp{3cm}|p{4.05cm}|p{4cm}}
(a)& finite subgroup $G$ of $\SL_2(\C)$ & 
(nontrivial) irreducible re\-pre\-sen\-tations &
decompositions of tensor products\\
\hline
(b)& simple Lie algebra of type $ADE$ & simple roots & 
(extended) Cartan matrix \\
\hline
(c)& minimal resolution $X\to\C^2/G$ &
irreducible components of the exceptional set,
or a basis of $H_2(X,\mathbb Z)$ &
intersection matrix
\end{tabular}
\end{center}

Gonzales-Sprinberg and Verdier~\cite{GV} realized the correspondence
between (a) and (c) geometrically as follows.
Let us consider the diagram
\begin{equation*}
\begin{CD}
  X @<p<< \uZ @>q>> \C^n,
\end{CD}
\end{equation*}
where $\uZ\subset X\times\C^n$ is the universal subscheme and $p$ and
$q$ are the projections to the first and second factors.
Let us define the {\it tautological bundle\/} $\tautR$ by
\begin{equation}\label{eq:deftaut}
  \tautR \defeq p_*\shfO_\uZ.
\end{equation}
Since $\uZ$ has a $G$-action, each fiber of $\tautR$ has a structure of
a $G$-module. By (3) in the definition of $X$, it is the regular
representation. We decompose into irreducibles:
\begin{equation}\label{eq:bundec}
   \tautR = \bigoplus_k \tautR_k \otimes \rho_k.
\end{equation}
Then Gonzales-Sprinberg and Verdier observed that
\begin{enumerate}
\item $\{\tautR_k\}_{k=0}^r$ gives a basis of the Grothendieck group
$K(X)$ of algebraic vector bundles over $X$,
\item $\{ c_1(\tautR_k) \}_{k\neq 0}$ is the dual basis of 
$\{ [ C_k ] \}$.
\end{enumerate}

Based on their results and also calculations by Vafa~et~al.\ related
to the mirror symmetry, Reid conjectured an existence of a similar
correspondence between (a) and (c), when $X$ is a crepant resolution of
$\C^n/G$ for $G\subset\SL_n(\C)$. (See \cite{Maizuru} for the history,
earlier results, and concrete examples of the McKay correspondence. We
do not reproduce them here.)

In this paper, we give more precise formulation of the conjecture when
the crepant resolution is given as $X$ defined above, and verify this
new conjecture when $G$ is a subgroup of $\SL_2(\C)$, or an abelian
subgroup of $\SL_3(\C)$.

Our new point is to consider the Grothendieck group of bounded
complexes of algebraic vector bundles with supports contained in
$\pi^{-1}(0)$, denoted by $K^c(X)$. There are natural elements $S_k$
of $K^c(X)$ which are also indexed by irreducible representations as
follows. The multiplication of the coordinate functions
$(x_1,\dots,x_n)$ on $\C^n$ induces the $G$-equivariant homomorphism
(called {\it tautological homomorphism\/})
\begin{equation*}
B \colon \tautR \to Q\otimes \tautR,
\end{equation*}
where $Q$ is the $n$-dimensional representation given by the
inclusion $G\subset\SL_n(\C)$. It gives rise the complex
\begin{equation}\label{eq:cpxnil}
\begin{CD}
\tautR @>d_n>> Q\otimes \tautR @>d_{n-1}>> \cdots @>d_2>> 
\bigwedge^{n-1}Q\otimes\tautR @>d_1>> \bigwedge^nQ\otimes\tautR = \tautR,
\end{CD}
\end{equation}
where $d_k(\eta) = B\wedge\eta$. This complex decomposes according to
\eqref{eq:bundec} and consider its transpose
\begin{equation}\label{eq:Sk}
S_k \colon \tautR_k^\vee \longrightarrow
\bigoplus_l a^{(n-1)}_{kl} \tautR_l^\vee 
\longrightarrow \cdots \longrightarrow
\bigoplus_l a^{(1)}_{kl} \tautR_l^\vee 
\longrightarrow \tautR_k^\vee,
\end{equation}
where the coefficients $a^{(i)}_{kl}$ is determined by 
\begin{equation}\label{eq:tensor}
{\textstyle\bigwedge^i Q}\otimes\rho_l = \bigoplus_k a^{(i)}_{kl} \rho_k.
\end{equation}
This is a generalization of the tensor product considered by
McKay.

We will show that $\{ \tautR_k \}_{k=0}^r$ and $\{ S_k \}_{k=0}^r$
form dual bases of $K(X)$ and $K^c(X)$ under the above assumption on
$G$. And we conjecture it holds for arbitary $G\subset\SL_3(\C)$.
Then from this approach, it becomes clear that the intersection
product among $S_k$'s are related to the decomposition of the tensor
product (see \corref{cor:int} for more precise statement).
Thus, our approach gives a `natural' explanation of the reason why the
decomposition of the tensor product is identified with the
intersection products in dimension $2$. As far as we know, known
proofs of this identification in dimension $2$ used case-by-case
analysis except those given in \cite[Appendix]{KrNa} and
\cite{Na-quiver}. Our proof is more natural and can be generalized to
higher dimensions.

The most essential ingredient in the proof of our main theorem is a
construction of a certain complex (see \eqref{eq:kosvec}).
We conjecture that it gives rise a resolution of the diagonal in
$X\times X$ for any $G\subset\SL_3(\C)$, and prove it when $G$ is
abelian. This complex is an
analogue of the Koszul complex on $\C^n$, and consists of vector
bundles of forms
\begin{equation*}
  \bigoplus_a p_1^* E_a\otimes p_2^* F_a,
\end{equation*}
where $p_1$, $p_2\colon X\times X\to X$ are projections into the first
and second factor. It is a higher-dimensional generalization of the
complex introduced by Kronheimer in his joint work with the second
author in $2$-dimensional case~\cite{KrNa}. The above conjecture was
proved there.

If $X$ would be compact, a standard argument (cf.\ 
\cite[Theorem~2.1]{ES}) shows that the Grothendieck group $K(X)$ of
vector bundles on $X$ is generated by $E_a$, either $F_a$.  However
$X$ is not compact, so the argument does not apply to our
situation. To overcome this difficulty, we modify the complex to
\begin{equation*}
\bigoplus_k p_1^* S_k \otimes p_2^*\tautR_k,
\end{equation*}
where $\tautR_k$ and $S_k$ are as above. The original complex and this
new complex are connected by a homotopy and defines the same element
in the Grothendieck group. This will lead us to our main theorem
(\thmref{thm:main}). The usage of both $K(X)$ and $K^c(X)$ are quite
essential in the argument.

Let us comment on our assumption on $G$. As we explained, we need this
assumption to show the exactness of the analogue of the Koszul
complex, more precisely, the condition~\eqref{eq:mainassum}. When $G$
is abelian, we can use the torus action so that we need to check
\eqref{eq:mainassum} for very specific ideals. The idea to use the
torus action is due to Nakamura~\cite{Iku} who used it to prove the
smoothness of $X$.

By the way, the correspondence between (b) and (c) was further developed
by the second author~\cite{Na-Duke,Na-affine}. He constructed irreducible
integrable representations of the affine Lie algebra on the homology
group of moduli spaces of instantons on $X$. The corresponding result
in dimension $3$ remains untouched in this paper.

The paper is organized as follows.
In \secref{sec:hilb}, we prepare some results on Hilbert schemes of
points on $\C^n$ and the fixed point component $X$ of the Hilbert
scheme.  In \secref{sec:quiver}, we identify $X$ with the moduli space
of a certain quiver. In \secref{sec:koszul}, we define the complex on
$X\times X$, and show that it gives a resolution of the diagonal
$\Delta$ under the condition~\eqref{eq:mainassum}. In
\secref{sec:McKay}, we state our main results on McKay correspondence
which is a correspondence between the representation ring $R(G)$ and
the Grothendieck group $K(X)$. We also study $K^c(X)$, the
Grothendieck group of bounded complexes of algebraic vector bundles
over $X$.  In \secref{sec:2d}, we study 2 dimensional case for the
argument in previous section.  In \secref{sec:messy}, we check the
condition holds for abelian subgroups $G\subset\SL_3(\C)$ and complete
the proof for our $3$ dimensional McKay correspondence for abelian
groups.

\subsection*{Acknowledgments}
The authors would like to express their deep gratitude to Professor
I.~Nakamura for discussion during the preparation of the paper.

Parts of this paper were written when the first author stayed at
University of Mannheim, and the second author stayed at University of
Bielefeld with support by SFB~343. Both authors express hearty thanks
to both universities for their hospitality.

\section{Fixed points in Hilbert schemes}\label{sec:hilb}

In this section, we prepare some preliminary results on Hilbert
schemes of points on $\C^n$ and the variety $X$ defined in the
introduction.

For a positive integer $N$, let $\Hilb{\C^n}{N}$ be the Hilbert scheme
parametrising $0$-dimensional subschemes of length $N$ (see
\cite{Lecture} for survey on Hilbert schemes of points). In this
paper, we shall confuse subschemes and the corresponding ideal of the
ring $\C[x_1, \dots, x_n]$. A point in $\Hilb{\C^n}{N}$ is either a
zero dimensional subscheme $Z\subset\C^n$ or an ideal $I\subset
\C[x_1,\dots,x_n]$.

Let $G$ be a finite subgroup of $\SL_n(\C)$. We consider two types of
`quotients' of $\C^n$ divided by $G$. The first one is the usual
set-theoretical quotient $\C^n/G$. It is a subvariety of the $N$th ($N
= \# G$) symmetric product $S^N(\C^n)= (\C^n)^N/{\mathfrak S}_N$
by the embedding
\begin{equation*}
  \C^n/G \ni G x \longmapsto \sum_{g\in G} [gx] \in S^N(\C^n),
\end{equation*}
where a point in the symmetric product is denoted by a formal sum of
points, as usual. The symmetric product is the Chow scheme of $N$
points in $\C^n$ parametrising effective $0$-cycles. Hence $\C^n/G$ is
the {\it Chow quotient\/} in the sense of \cite{Kap}. It is an
irreducible component of the fixed point of the induced $G$-action on
$S^N(\C^n)$.

Another quotient is the {\it Hilbert quotient\/} which is obtained by
replacing the symmetric product by the Hilbert scheme as follows:
Consider the induced $G$-action on $\Hilb{\C^n}{N}$. If
$Z\in\Hilb{\C^n}{N}$ is a fixed point, then $H^0(\shfO_Z)$ is a
$G$-module. For example, if $Z$ is a single $G$-orbit consisting of
pairwise distinct $N$ points, $H^0(\shfO_Z)$ is the regular
representation of $G$.
As in \secref{sec:intro}, let $X$ be the variety parametrising
$Z\in(\Hilb{\C^n}{N})^G$ such that $H^0(\shfO_Z)$ is the regular
representation of $G$. 
The $G$-module structure is constant on each connected component of
the fixed point set $(\Hilb{\C^n}{N})^G$. Thus $X$ is a union of
components. A priori, it may consists of several irreducible
components, but $X$ is irreducible in many cases as we will see later.

We have the {\it Hilbert-Chow\/} morphism $\pi$ from
$\Hilb{\C^n}{N}$ to $S^N(\C^n)$ defined by
\begin{equation*}
  \pi\colon \Hilb{\C^n}{N}\ni Z \longmapsto \sum_{x\in \C^n}
  \operatorname{length}(Z_x)[x]\in S^N(\C^n).
\end{equation*}
Take a point $Z$ in $X$ and consider $\pi(Z)$. Since $Z$ is invariant
under $G$, its support consists of a union of $G$-orbits. However,
since constant functions on each orbit form the trivial representation
contained in $H^0(\shfO_Z)$, we only have single $G$-orbit by the
assumption that $H^0(\shfO_Z)$ is the regular representation. This
implies that $\pi(Z)$ is of the form
\begin{equation*}
  \sum_{g\in G} [gx] = \sum_{y\in Gx} \#(G/G_x) [y],
\end{equation*}
for some $x\in\C^n$. Here $Gx$ denote the $G$-orbit through $x$ and
$G_x$ is the stabilizer of $x$ in $G$. Hence $\pi$ maps $X$ to
$\C^n/G$. We use the same notation $\pi$ for the restriction of the
map to $X$ for brevity.

The nonsingular locus $(\C^n/G)^{\operatorname{reg}}$ of $\C^n/G$
consists of those orbits $Gx$ with $G_x$ trivial. Since the map $\pi$
is an isomorphism on $\pi^{-1}((\C^n/G)^{\operatorname{reg}})$,
$\pi\colon X\to \C^n/G$ is a resolution of singularities provided
$X$ is nonsingular of dimension $n$ and connected.

\begin{Remark}
(1) The terminologies, Chow quotients and Hilbert quotients, were
introduced by Kapranov~\cite{Kap}.

(2) The definition of $X$ given in \cite{GK,IN,Iku} is slightly
different from above. In those papers, $X$ is defined as the
irreducible components of $(\Hilb{\C^n}{N})^G$ containing $G$-orbits
of cardinality $N$. In dimension $2$, it is known that the above $X$
is smooth and connected, and hence two definitions are same (see
\cite[4.4]{Lecture}). We prove the connectedness of $X$ when $G$
is an abelian subgroup of $\SL_3(\C)$ later. Thus the definition is also
the same in this case.
\end{Remark}

\section{Representations of a quiver}\label{sec:quiver}

In this section, we identify the subvariety $X$ of the
Hilbert scheme with the moduli space of a certain quiver. This
identification shows that $X$ is a special case of the variety
considered by Kronheimer~\cite{Kr} (in dimension $2$) and
Sardo-Infirri~\cite{Sacha} (in general). We hope that this section
will be helpful for the reader to notice that the complex constructed
in the next section is a natural generalization of the complex
introduced by Kronheimer-Nakajima~\cite{KrNa}.

Let $A = \C[x_1,x_2,\dots, x_n]$ be the coordinate ring of $\C^n$.
Take $I\in X$, or more generally an ideal corresponding to a
zero-dimensional subscheme of $\C^n$ of length $N$. Then $A/I$ is an
$N$-dimensional vector space. If $Z$ is the corresponding subscheme,
we have $A/I = H^0(\shfO_{Z}$). The multiplications of the coordinate
functions $(x_1,x_2,\dots,x_n)$ induce a map
$B = (B_1,B_2,\dots, B_n) \colon A/I \to \C^n\otimes A/I$ by
\begin{equation*}
  B_\alpha(f\bmod I) \defeq x_\alpha f\bmod I, \qquad
  (\alpha = 1,2,\dots, n).
\end{equation*}
It satisfies
\begin{equation*}
\Hom(A/I, {\textstyle \bigwedge^2}\C^n\otimes A/I)\ni
[B\wedge B] =
\sum_{\alpha < \beta} [B_\alpha, B_\beta] dx_\alpha\wedge dx_\beta = 0.
\end{equation*}
Let us define $i\colon \C \to A/I$ by
$i(\lambda) = \lambda\bmod I$. Then $i(1)$ is a cyclic vector with
respect to $B_\alpha$'s, that is there is no proper subspace $S\subsetneq
A/I$ which contains $i(1)$ and is invariant under all $B_\alpha$'s.

Conversely if we have an $N$-dimensional vector space $R$ and
homomorphisms $B\colon R \to \C^n\otimes R$, $i\colon \C\to R$ such
that $[B\wedge B] = 0$, and $i(1)$ is a cyclic vector with respect to
$B_\alpha$'s, we can define an ideal $I$ by
\begin{equation*}
  I \defeq \{ f(x_1,\dots, x_n) \in A \mid f(B_1,\dots,B_n)i(1) = 0\}.
\end{equation*}
Then $I$ defines a $0$-dimensional subscheme of $\C^n$ of length $N$.

Now suppose $I$ is a point in $X$, i.e.\  a) it is invariant under the
action of $G$ and b) $A/I$ is the regular representations of $G$.
Then the above homomorphisms $B\colon
A/I\to Q\otimes A/I$, $i\colon \C\to A/I$ are $G$-equivariant,
where $Q$ is a $G$-module defined by the inclusion
$G\subset\SL_n(\C)$, and $\C$ is the trivial $G$-module. Hence we get
\begin{Proposition}\label{prop:quiver}
Let $R$ be the regular representation of $G$. Then
there exists a bijection between $X$ and
the quotient space of the homomorphisms
$B\in \Hom_G(R, Q\otimes R)$, $i\in \Hom_G(\C, R)$ satisfying
\begin{align}
  & [B\wedge B] = 0, \label{eq:ADHM}\\
  & \text{$i(1)$ is a cyclic vector
  with respect to $B_\alpha$'s}\label{eq:cyclic}
\end{align}
by the action of $\GL_G(R)$, the group of $G$-equivariant
automorphisms of $R$.
\end{Proposition}

Let us rewrite the above quotient space as a moduli space of
representations of a certain quiver. Let $\rho_0$, \dots, $\rho_r$ be
the isomorphism classes of irreducible representations of $G$, with
$\rho_0$ be the trivial representation. Then the regular
representation $R$ decomposes as
\begin{equation}\label{eq:irrdec}
  R = \bigoplus R_k\otimes \rho_k,
\end{equation}
where $R_k = \Hom_G(\rho_k, R)$.
Then we have
\begin{equation*}
  \Hom_G(R, Q\otimes R) = \bigoplus \Hom(R_k, R_l)\otimes
  \Hom_G(\rho_k, Q\otimes \rho_l) =
  \bigoplus a_{kl}\Hom(R_k, R_l),
\end{equation*}
where $a_{kl} = a^{(1)}_{kl}$ is given in \eqref{eq:tensor}.
Similarly we have
\begin{equation*}
  i \in \Hom_G(\C,R) = \Hom(\C, R_0).
\end{equation*}
The group $\GL_G(R)$ of $G$-equivariant automorphisms of $R$ can be
rewritten as
\begin{equation*}
  \GL_G(R) = \prod_k \GL(R_k).
\end{equation*}

Hence $X$ can be described as
\begin{equation}\label{eq:decdes}
  \{ (B,i)\in \bigoplus a_{kl}\Hom(R_k, R_l)\oplus\Hom(\C, R_0)\mid
\eqref{eq:ADHM},\eqref{eq:cyclic} \} / \prod_k \GL(R_k).
\end{equation}
This description depends only on $a_{kl}$ and $\dim R_k = \dim \rho_k$.

\begin{figure}[htbp]
\begin{center}
{\includegraphics{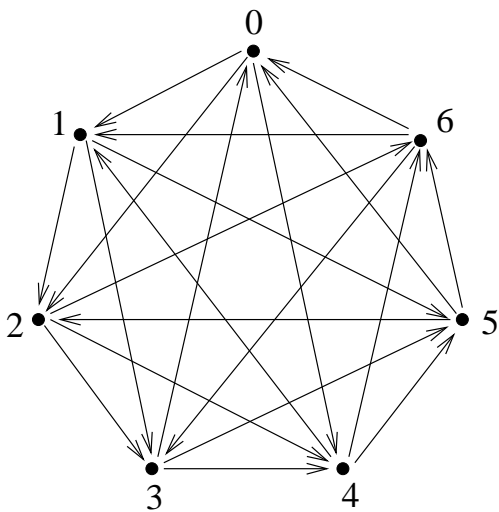}}
\caption{McKay quiver for $G = \langle\
\operatorname{diag}(\varepsilon,\varepsilon^2,\varepsilon^4)\ \rangle$
($\varepsilon = \exp(2\pi i/7)$)}
\label{fig:McKay}
\end{center}
\end{figure}

The {\it McKay quiver\/} is the quiver whose vertices are irreducible
representation with $a_{kl}$ arrows (possibly $0$) from the vertex $k$
to the vertex $l$. The above space is the framed moduli space of
representation of the McKay quiver with the relation corresponding to
\eqref{eq:ADHM}. (cf.\ \cite{Na-quiver}) (See Figure~\ref{fig:McKay}.)

When $n=2$, the description in \propref{prop:quiver} 
is essentially same as Kronheimer's construction of
ALE spaces~\cite{Kr}.
There are minor differences: First his space depends on a parameter
$\zeta$. Our space corresponds his space with a special choice of
$\zeta$. Second, we have an extra vector $i$ and take
quotient by $\GL_G(R)$, while Kronheimer had no $i$ and took quotients by
$\GL_G(R)/\text{scalar}$.  However, we can always normalize as $i =
1$ by the action of $\GL(R_0) = \C^\ast$. Hence our quotient is
isomorphic to the quotient space of $B$'s by the action of
$\prod_{k\ne 0} \GL(R_k) \cong \GL_G(R)/\text{scalar}$. Hence our
description is same as \cite{Kr}.

Kronheimer's construction was generalized to higher dimensions by
Sardo-Infirri \cite{Sacha}. Thus our description coincides with his
with a particular parameter.

\section{A Koszul complex}\label{sec:koszul}

In this section we construct a resolution of the diagonal
$\Delta$ in $X\times X$ following \cite[3.6]{KrNa}.

Take $I_1$, $I_2\in X$, and consider the corresponding $G$-equivariant
homomorphisms $A/I_1\to Q\otimes A/I_1$, $A/I_2\to Q\otimes A/I_2$
as in the previous section. Let us denote them by $B^1$, $B^2$.
Then consider the following complex of vector spaces:
\begin{equation}\label{eq:koszul}
\begin{CD}
  E_n @>d_n>> E_{n-1} @>d_{n-1}>> \cdots @>d_2>> E_1 @>d_1>> E_0,
\end{CD}
\end{equation}
where
\begin{align*}
  & E_k \defeq \Hom_G(A/I_1, \textstyle{\bigwedge^{n-k}}Q\otimes A/I_2), \\
  & d_k (\eta) \defeq B^2\wedge \eta - \eta\wedge B^1\quad
  \text{for $\eta\in E_k$}.
\end{align*}
The equality $d_{k-1}\circ d_k = 0$ follows from the equation
$[B^a\wedge B^a] = 0$ ($a=1,2$).

\begin{Lemma}[cf.\ \protect{\cite[3.8,3.9]{KrNa}}]\label{lem:H0}
\rom{(1)} When $I_1 \neq I_2$, $d_n$ is injective and $d_1$ is
surjective.

\rom{(2)} When $I_1 = I_2$, the kernel of $d_n$ is the one-dimensional
subspace of scalar endomorphisms and the image of $d_1$ is the
codimension-one subspace of trace-free endomorphisms.
\end{Lemma}

\begin{proof}
First notice that $\bigwedge^n Q$ is the trivial $G$-module by the
assumption $G\subset\SL_n(\C)$. Then we have $\bigwedge^{n-1} Q \cong
Q^*$, and the statements for $d_1$ can be proved dually by the same
arguments as those for $d_n$. So we only give the proof on the
statements for $d_n$.

Suppose $\eta$ is an element in the kernel of $d_n$. We have
\begin{equation}\label{eq:equiv}
  \eta B^1 = B^2 \eta.
\end{equation}
This equation implies that the image of $\eta$ is invariant under
$B^2$.

Since $\eta$ is $G$-equivariant, the $G$-fixed parts are preserved
under $\eta$. Since $A/I_1$ and $A/I_2$ are regular representations,
the $G$-fixed parts consist of constant multiples of $1\bmod I_1$ and
$1\bmod I_2$ respectively. Hence, we have $\eta(1\bmod I_1) = \lambda
(1\bmod I_2)$ for some constant $\lambda$. Then \eqref{eq:equiv}
implies
\begin{equation*}
  \eta(x_1^{k_1}\cdots x_n^{k_n}\bmod I_1) =
    \lambda x_1^{k_1}\cdots x_n^{k_n}\bmod I_2
\end{equation*}
for any $k_1$,\dots, $k_n\in\mathbb Z_{\ge 0}$. If $\lambda = 0$, then
$\eta = 0$. If $\lambda \neq 0$, then $\eta$ is surjective. Since
$A/I_1$ and $A/I_2$ has the same dimension, this implies $I_1 = I_2$
and $\eta$ is a scalar endomorphism.
\end{proof}

\begin{Lemma}[cf.\ \protect{\cite[3.10]{KrNa}}]\label{lem:H1}
\rom{(1)} When $I_1 = I_2$, the homology group
$\Ker d_{n-1}/\Ima d_n$ is isomorphic to $\Hom_A(I_1, A/I_2)^G$, the
$G$-fixed part of $\Hom_A(I_1, A/I_2)$.

\rom{(2)} When $I_1\neq I_2$,
$\Ker d_{n-1}/\Ima d_n$ is isomorphic to $\Hom_A(I_1, A/I_2)^G /
\C \Xi$, where $\Xi$ is the composition
\begin{equation*}
   \Xi\colon I_1 \hookrightarrow A \to A/I_2.
\end{equation*}

Moreover, in either case, $\Ker d_1/\Ima d_2$ is isomorphic to dual
space of $\Ker d_{n-1}/\allowbreak\Ima d_n$ of the complex with $I_1$,
$I_2$ are exchanged.
\end{Lemma}

\begin{proof}
The statement for the duality between the degree $n-1$ and degree
$1$ follows exactly as in the proof of \lemref{lem:H0}.

Let $\Phi\in \Hom_A(I_1, A/I_2)^G$. We take an extension $\Psi\colon
A\to A/I_2$ as a $G$-equivariant homomorphism (not necessarily
an $A$-homomorphism). Then we define $\dot{B}
= (\dot{B}_1, \dots, \dot{B}_n)\in\Hom_G(A/I_1, Q\otimes A/I_2)$ by
\begin{equation*}
  \dot{B_\alpha}(f \bmod I_1) \defeq \Psi(x_\alpha f) -
  B^2_\alpha\Psi(f) \qquad\text{for $\alpha=1,\dots, n$, $f\in A$}.
\end{equation*}
When $f\in I_1$, the right hand side vanishes since $\Phi$ is an
$A$-homomorphism. Hence $\dot{B}$ is well-defined.

It is easy to check $\dot{B}\in \Ker d_{n-1}$.
Moreover, the ambiguity of the choice of
the extension $\Psi$ is compensated by the image of $d_n$.
Hence we have a homomorphism
\begin{equation}
  \Hom_A(I_1, A/I_2)^G \to \Ker d_{n-1}/\Ima d_n.
\label{eq:cohom}\end{equation}

Conversely, if we are given $\dot{B}\in\Ker d_{n-1}$, we define 
an $G$-equivariant homomorphism $\Psi\colon A\to A/I_2$ inductively by
\begin{equation*}
\left\{\begin{split}
  & \Psi(1) = 1\bmod I_2\\
  & \Psi(x_\alpha f) = B_\alpha^2\Psi(f) + \dot{B}_\alpha(f\bmod I_1)
  \qquad\text{for $\alpha=1,\dots, n$, $f\in A$}.
\end{split}\right.
\end{equation*}
This is well-defined thanks to the assumption $d_{n-1}\dot{B} =
0$. Moreover, the restriction $\Phi = \Psi|_{I_1}$ is $A$-linear by
the second equation. This argument shows that the map \eqref{eq:cohom} 
is surjective.

Now suppose $\Phi$ lies in the kernel of \eqref{eq:cohom}. Then we can 
take the extension $\Psi\colon A\to A/I_2$ of $\Phi$ so that
\begin{equation*}
  \Psi(x_\alpha f) = B_\alpha^2 \Psi(f)\qquad\text{for
  $\alpha=1,\dots, n$, $f\in A$}.
\end{equation*}
Since $\Psi$ is $G$-equivariant, we have $\Psi(1) = \lambda \bmod I_2$
for some constant $\lambda$. Then the above equation implies that
$\Psi(f) = \lambda f\bmod I_2$ for any $f\in A$. Hence the restriction
$\Phi = \Psi|_{I_1}$ is $\lambda\Xi$. If $I_1 = I_2$, then $\Xi$
is zero, and hence we have the assertion.
\end{proof}

Note that $\Hom_A(I, A/I)^G$ is the Zariski tangent space of $X$ at
$I$. In particular, $X$ is nonsingular if and only if
$\Hom_A(I, A/I)^G$ has a constant dimension independent of $I$ on each
irreducible component.

\begin{Corollary}\label{cor:bundle}
Suppose $X$ is nonsingular. 
Then the bijection given in \propref{prop:quiver} is an isomorphism.
Moreover, the variety $P$ consisting of pairs
$(B,i)\in \Hom_G(R, Q\otimes R)\times \Hom_G(\C, R)$ \rom($R$ is the
regular representation of $G$ as before\rom) satisfying both
\eqref{eq:ADHM}, \eqref{eq:cyclic}
is nonsingular, and the quotient map $P \to P/\GL_G(R) = X$ is a
$\GL_G(R)$-principal bundle.
\end{Corollary}

\begin{proof}
Consider the complex \eqref{eq:koszul} for $I_1 = I_2$ and $R =
A/I_1$. Then $d_{n-1}$ is nothing but the differential of the map
$B\mapsto [B\wedge B]$.

The assumption of the smoothness of $X$ implies that $\Hom_A(I_1,
A/I_1)^G$, the Zariski tangent space at $I_1$, has a constant
dimension independent of $I_1$. By \lemref{lem:H0} and
\lemref{lem:H1}, the kernel of $d_1$ also has a constant
dimension. Hence the variety $P$ is nonsingular.

The action of $\GL_G(R)$ on $P$ is free. For, if $g\in\GL_G(R)$
stabilizes $B$ and $i$, then $\Ker(g-1)$ contains $i(1)$ and is
invariant under $B$. Hence the cyclic vector condition implies 
$\Ker(g-1) = R$.

Moreover, $P/\GL_G(R)$ is a geometric invariant theory quotient of $P$
by $G$ if we introduce the polarization as in the case of quiver
varieties \cite{Na-affine}. The cyclic vector condition is the
stability in the geometric invariant theory, and hence the
quotient map is a $\GL_G(R)$-principal bundle.

Now the map given by the \propref{prop:quiver} is differentiable and
respects tangent spaces. Hence it is an isomorphism.
\end{proof}

Assume that $X$ is nonsingular of dimension $n$. Then if we vary $I_1$
in $X$, $A/I_1$ forms a holomorphic vector bundle over $X$. In fact,
it is identified with the associated vector bundle $P\times_{\GL_G(R)}
R$. Following \cite{KrNa}, we denote it by $\tautR$ and call {\it
tautological vector bundle}. Fibers of $\tautR$ have structures of
$G$-modules which are isomorphic to the regular representation. The
homomorphism $B\in \Hom_G(A/I, Q\otimes A/I)$ defines a
$G$-equivariant holomorphic vector bundle homomorphism $\tautR \to
Q\otimes\tautR$. This is called {\it tautological endomorphisms\/} in
\cite{KrNa}.

Since $A/I$ is the $0$-th cohomology of the structure sheaf of
subscheme corresponding to $I$, we can identify the tautological
bundle $\tautR$ with $p_*\shfO_\uZ$ where $\uZ \subset X\times \C^3$
is the universal subscheme and $p\colon \uZ\to X$ is the first
projection.  Thus we arrive at the definition~\eqref{eq:deftaut}.

As explained in the introduction, the decomposition of the regular
representation into irreducible representations~\eqref{eq:irrdec}
induces the decomposition of the tautological vector bundle as
\eqref{eq:bundec}. In other words,
\begin{equation*}
   \tautR_k 
    = \Hom_G(\rho_k, p_*\shfO_\uZ).
\end{equation*}

For brevity, we identify the vector bundle $\tautR$ with the
sheaf of germs of its sections. Then the complex \eqref{eq:koszul}
induces the following complex of sheaves on $X\times X$:
\begin{equation}\label{eq:kosvec}
\begin{CD}
  0 @>>> F_n @>d_n>>\cdots @>d_2>> F_1 @>d_1>> F_0
  @>t>> \shfO_{\Delta} @>>> 0,
\end{CD}
\end{equation}
where
\begin{align*}
  & F_k \defeq \Hom_G(p_1^* \tautR,
  \textstyle{\bigwedge^{n-k}}Q\otimes p_2^*\tautR), \\
  & d_k (\eta) \defeq B^2\wedge \eta - \eta\wedge B^1\quad
  \text{for $\eta\in F_k$}, \\
  & t(\eta) = \tr(\eta|_\Delta) \quad
  \text{for $\eta\in F_0$}.
\end{align*}
Here $p_a$ is the projection to the $a$th factor of $X\times X$. Note
that $F_0 = \Hom_G(p_1^*\tautR, p_2^*\tautR)$ and hence the trace
make sense on the diagonal.

When $n=2$, the Hilbert scheme of points is nonsingular by
Fogarty~\cite{Fog}. Hence the $G$-fixed component $X$ is also
nonsingular. Thus the assumption of \corref{cor:bundle} is
met. Counting dimensions and using \lemref{lem:H0}, we deduce that
$\Ker d_1/\Ima d_2 = 0$ for $I_1\neq I_2$. This and an additional
argument shows that the complex~\eqref{eq:kosvec} is exact (see
\cite[\S3]{KrNa} or below).
Moreover, the last assertion in \lemref{lem:H1} says that $\Ker
d_1/\Ima d_2$ for $I_1 = I_2$, which is the tangent space of $X$ at
$I_1$, is isomorphic to its dual space. This isomorphism is given by
the natural holomorphic symplectic form on $X$.

We assume $n = 3$ hereafter, and consider the following condition:
\begin{equation}\label{eq:mainassum}
 \dim \Hom_A(I_1, A/I_2)^G =
 \begin{cases}
 3, & \text{when $I_1 = I_2$}\\
 1, & \text{when $I_1 \neq I_2$}.
 \end{cases}
\end{equation}
When $I_1$, $I_2$ are ideals given by distinct points, the above
holds.  The condition for $I_1 = I_2$ is equivalent to saying that $X$
is nonsingular of dimension $3$. And note that similar condition holds
for $n=2$ by the above discussion. Those show that the above seems reasonable. And we show the above holds when $G$ is abelian in
\secref{sec:messy}.

\begin{Theorem}[cf.\ \protect{\cite[3.6]{KrNa}}]\label{thm:exact}
Under the assumption~\eqref{eq:mainassum}, the
complex~\eqref{eq:kosvec} is exact.
\end{Theorem}

\begin{proof}
By the assumption, the complex is exact outside the diagonal $\Delta$.
The exactness of \eqref{eq:kosvec} in degree $0$ can be shown exactly
as in \cite[3.6]{KrNa}. So we omit the argument.

For the proof of the exactness in degrees other than $0$, we use the
criterion of Buchsbaum-Eisenbud (see e.g., \cite[20.9]{Eis}). We need
to check (a) $\rank d_k + \rank d_{k+1} = \rank F_k$, and (b) the
determinantal ideal of the differential $d_k$ has depth at least
$k$. The first condition holds since the complex is exact on a
nonempty open subset. Since the diagonal has codimension $3$ in
$X\times X$, the determinantal ideal of the differential $d_k$ has
depth at least $3$. Hence the second condition (in fact, a stronger
condition) also holds.
\end{proof}

The following will not used in the other part of this paper, but
illustrates the relation between the smoothness and the exactness of
the complex~\eqref{eq:kosvec}.

\begin{Proposition}
If $X$ is nonsingular, the complex
\eqref{eq:kosvec} is exact on the diagonal $\Delta\subset X\times X$.
\label{prop:exact}\end{Proposition}

\begin{proof}
Fix a point $x_0\in X$ and consider the complex \eqref{eq:kosvec} at the
point $(x_0, x_0)$:
\begin{equation}\label{eq:compt}
\begin{CD}
  (F_3)_{(x_0,x_0)} @>(d_3)_{(x_0,x_0)}>> (F_2)_{(x_0,x_0)}
  @>(d_2)_{(x_0,x_0)}>> (F_1)_{(x_0,x_0)} @>(d_1)_{(x_0,x_0)} >>
  (F_0)_{(x_0,x_0)}.
\end{CD}
\end{equation}
By Lemmas~\ref{lem:H0},\ref{lem:H1} together with the smoothness
assumption of $X$, the homology groups of this complex are
$\C(\operatorname{Id})$ in degree $0$ and $3$, the tangent space
$T_{x_0}X$ in degree $2$, and the cotangent space $T^*_{x_0}X$ in
degree $1$.
We define a trivial vector bundle $H_i$ over the tangent space
$T_{x_0}X\times T_{x_0}X$ where the fiber is the $i$th homology
groups of \eqref{eq:compt}.
 
Choosing connections on vector bundles in \eqref{eq:kosvec}, we
consider the derivative $D d_k$ at $x_0$. Differentiating
$d_{k-1}\circ d_k = 0$, we check that $D d_k$'s induce homomorphisms
between homology groups of \eqref{eq:compt}. Moreover, they are
independent of the choice of the connections. Let us think the
homomorphisms as vector bundles homomorphisms between vector bundles
$H^k$ by setting the value at $(v,w)\in T_{x_0}X\times T_{x_0}X$ as
the derivative $D d_k$ in the direction $(v,w)$. We simply write $D
d_k$ for these vector bundle homomorphisms.
Thus we get
\begin{equation}\label{eq:comtan}
\begin{CD}
  0 @>>> H_3 @>Dd_3>> H_2 @>Dd_2>>
  H_1 @>Dd_1>> H_0
\end{CD}
\end{equation}
Differentiating $d_{k-1}\circ d_k = 0$ twice, one checks that this
forms a complex.

Calculating the derivative of each $d_k$, we find that
\eqref{eq:comtan} is a part of Koszul complex on $T_{x_0} X \times
T_{x_0}X \cong \C^3\times \C^3$. Namely, if we add the evaluation
homomorphism $H_0 \to \shfO_\Delta$ to the end of \eqref{eq:comtan},
the complex is the Koszul complex. Now we check the exactness of the
original complex~\eqref{eq:kosvec} at $(x_0,x_0)$, except possibly in
degree $0$. For this we again use the criterion of
Buchsbaum-Eisenbud~\cite{Eis}. From the above, the ranks of $d_k$ are
\begin{gather*}
\rank d_1 = \rank (d_1)_{(x_0,x_0)} + 1, \quad
\rank d_2 = \rank (d_2)_{(x_0,x_0)} + 2, \\
\rank d_3 = \rank (d_3)_{(x_0,x_0)} + 1,
\end{gather*}
where $(d_k)_{(x_0,x_0)}$ is what $d_k$ induces on the fibers of $F_k$
over $(x_0,x_0)$. Combining with the above observation on the
homology groups of \eqref{eq:compt}, we get $\rank d_k + \rank
d_{k+1} = \rank F_k$. Moreover, the determinantal ideal of each
differential has depth $3$. These imply the exactness of the complex
as in \thmref{thm:exact}.
\end{proof}

\section{McKay correspondence for the $K$-theory}\label{sec:McKay}

Given a finite group $G$ acting on a variety $Y$, we denote by
$K_G(Y)$ the Grothendieck group of $G$-equivariant coherent
$\shfO_Y$-sheaves over $Y$. When $G$ is the trivial group $\{1\}$, we
simply write $K(Y)$ for $K_{\{1\}}(Y)$. This is the ordinary
$K$-group. We denote by $[S]$ the class represented by a
$G$-equivariant sheaf $S$ on $Y$. But we may drop the bracket when
there is no ambiguity. If $Y$ is nonsingular, $K_G(Y)$ is isomorphic
to the Grothendieck group of $G$-equivariant vector bundles (see e.g.,
\cite[Chapter 5]{CG}). If $f\colon Y\to Y'$ is a $G$-equivariant
proper morphism, we can define a push-forward $f_*\colon K_G(Y) \to
K_G(Y')$ by \( f_*([S]) = \sum_i (-1)^i [R^if_*(S)], \) where
$R^if_*(S)$ is the $i$th higher direct image sheaf. If $f\colon Y\to
Y'$ is a $G$-equivariant morphism and $Y'$ is smooth, we define the
pull-back $f^*\colon K_G(Y')\to K_G(Y)$ as follows: Since $Y'$ is
nonsingular, it is enough to define the pull-back for classes
represented by $G$-equivariant {\it vector bundles}. If $E$ is a
$G$-equivariant vector bundle on $Y'$, then its pull-back $f^*(E)$ is
also a $G$-equivariant vector bundle over $Y$. Hence we define
$f^*([E]) = [f^*(E)]$. We will never use the pull-back homomorphism
from singular varieties.

Let us consider a subvariety $\pi^{-1}(0)$ and $K^c(X)$ the
Grothendieck group of bounded complexes of algebraic vector bundles
over $X$ which are exact outside $\pi^{-1}(0)$. (See \cite{BFM} for the
definition and results used below.) This is isomorphic to the
Grothendieck group of coherent sheaves on $\pi^{-1}(0)$, where the
isomorphism is given by taking the alternating sum of the homology of
the complex:
\begin{equation*}
  [ E_\bullet ] = [ E_n \to E_{n-1} \to \cdots \to E_1 \to E_0 ]
  \longmapsto \sum_{i=0}^n (-1)^n [ H_i(E_\bullet) ].
\end{equation*}
The inverse is given by mapping a sheaf $S$ on $\pi^{-1}(0)$
to its finite resolution by locally free sheaves over $X$. When we
consider push-forward homomorphisms, we represent elements in $K^c(X)$
by sheaves on $\pi^{-1}(0)$. When we consider pull-back homomorphisms,
we represent elements by complexes.

There is a natural pairing between $K(X)$ and $K^c(X)$ given by
\begin{equation}\label{eq:pairing}
  K(X)\times K^c(X) \ni ([E], [S]) \mapsto
  P_*([E\otimes S])\in K(\operatorname{point})\cong \mathbb Z,
\end{equation}
where $E$ is a vector bundle on $X$ and $S$ is a sheaf on
$\pi^{-1}(0)$ and $P$ is the obvious projection of $\pi^{-1}(0)$ to a
point. Note that
\begin{aenume}
\item $E\otimes S$ is a tensor product of a vector bundle and a sheaf,
hence well-defined in the Grothendieck group,
\item $E\otimes S$ has support contained in $\pi^{-1}(0)$, hence
$P_*([E\otimes S])$ can be defined.
\end{aenume}

Let us consider the complex~\eqref{eq:cpxnil} in $n=3$:
\begin{equation*}
\begin{CD}
\tautR @>d_3>> Q\otimes \tautR @>d_2>> 
\bigwedge^{2}Q\otimes\tautR @>d_1>> \bigwedge^3Q\otimes\tautR = \tautR.
\end{CD}
\end{equation*}
It is a complex thanks to the equation $[B\wedge B] = 0$.

\begin{Lemma}\label{lem:support}
The complex \eqref{eq:cpxnil} is exact outside $\pi^{-1}(0)$.
\end{Lemma}

\begin{proof}
Take a coordinate system $(x_1, x_2, x_3)$ on $\C^3$ and write
$B = (B_1, B_2, B_3)$. Note that the support of the $0$-dimensional
subscheme corresponding to $[B_1, B_2, B_3, i]$ consists of
simultaneous eigenvalues of $B_1$, $B_2$, $B_3$.
Hence at least one of $B_\alpha$'s is invertible 
if $[B_1, B_2, B_3, i]$ is outside of $\pi^{-1}(0)$.
Say $B_1$ is invertible. Now it is clear that $d_3$ is injective and
$d_1$ is surjective.

Suppose $\eta = (\eta_1,\eta_2,\eta_3)$ is in the kernel of $d_2$,
that is
\begin{equation*}
  B_1\eta_2 = B_2\eta_1, \quad
  B_2\eta_3 = B_3\eta_2, \quad
  B_3\eta_1 = B_1\eta_3.
\end{equation*}
Setting $\xi = B_1^{-1}\eta_1$, we find
\begin{equation*}
  d_3 \xi = (\eta_1, B_2 B_1^{-1}\eta_1, B_3 B_1^{-1}\eta_1)
  = (\eta_1, \eta_2, \eta_3),
\end{equation*}
where we have used $[B_1, B_2] = [B_3, B_1] = 0$. This shows that 
$\Ker d_2 = \Ima d_3$.
The proof for $\Ker d_1 = \Ima d_2$ is same.
\end{proof}

We decompose the complex~\eqref{eq:cpxnil} according to
\eqref{eq:bundec} and denote its transpose by $S_k$:
\begin{equation*}
S_k \colon \tautR_k^\vee \longrightarrow
\bigoplus_l a^{(2)}_{kl} \tautR_l^\vee 
\longrightarrow
\bigoplus_l a^{(1)}_{kl} \tautR_l^\vee 
\longrightarrow \tautR_k^\vee.
\end{equation*}
By \lemref{lem:support}, $S_k$ defines an
element in $K^c(X)$.

Now we define the homomorphism from the representation ring $R(G)$ of
$G$ to $K(X)$ as follows.
Let us consider the diagram
\begin{equation*}
\begin{CD}
  X @<p<< \uZ @>q>> \C^n,
\end{CD}
\end{equation*}
where $\uZ\subset X\times\C^n$ is the universal subscheme and $p$ and
$q$ are the projections to the first and second factors. Note that the
group $G$ acts on $\uZ$ and $\C^n$ so that $q$ is $G$-equivariant. If
we let $G$ act on $X$ trivially, $p$ is also $G$-equivariant. By
\cite[5.4.21]{CG}, the representation ring $R(G)$ is isomorphic to
$K_G(\C^n)$ by sending the representation $\rho$ to
$\rho\otimes\shfO_{\C^n}$. We consider the composition of the
following homomorphisms in $K$-theory:
\begin{multline}\label{eq:Kcomp}
  R(G) \xrightarrow{\vee} R(G) \cong K_G(\C^n) \rightarrow \\
  \xrightarrow{q^*} K_G(\uZ) \xrightarrow{p_*} K_G(X)
  \cong R(G)\otimes_{\mathbb Z} K(X)
  \xrightarrow{\operatorname{Inv}\otimes\operatorname{id}} K(X),
\end{multline}
where $\vee$ is sending $V$ to its dual representation $V^\vee$, and
$\operatorname{Inv}\colon R(G)\to \mathbb Z$ is given by
$\operatorname{Inv}(V) = \dim V^G$.

The image of $\rho_k$ under the composition~\eqref{eq:Kcomp} is given
by
\begin{multline*}
  \left((\operatorname{Inv}\otimes\operatorname{id})
    \circ p_*\circ q^*\right)(\rho_k^\vee)
 = \Hom_G\left(\rho_0, (p_*\circ q^*)
        (\rho_k^\vee\otimes\shfO_{\C^3}))\right) \\
 = \Hom_G(\rho_k, p_*\shfO_\uZ) = \tautR_k.
\end{multline*}

The following is one of main results in this paper.
\begin{Theorem}\label{thm:main}
Suppose $G$ is a finite subgroup of $\SL_3(\C)$. Assume the
condition~\eqref{eq:mainassum} holds.

\rom{(1)} The composition $(\operatorname{Inv}\otimes\operatorname{id})
\circ p_*\circ q^*\circ\vee$ in \eqref{eq:Kcomp}, which maps the
irreducible representation $\rho_k$ to the tautological bundle
$\tautR_k$, gives an isomorphism between $R(G)$ and $K(X)$.

\rom{(2)} The support of the complex $S_k$ is contained in
$\pi^{-1}(0)$, and $\{\tautR_k\}$ and $\{S_k\}$ are dual bases for
$K(X)$ and $K^c(X)$, 
where $K^c(X)$ is the Grothendieck group of bounded complexes of vector
bundles with supports contained in $\pi^{-1}(0)$.
\end{Theorem}

The rest of this section is devoted to the proof of this theorem and to its
corollaries. So we assume the condition~\eqref{eq:mainassum}
throughout in this section. Hence the complex~\eqref{eq:kosvec} is
exact by \thmref{thm:exact}. The assumption~\eqref{eq:mainassum} will
be checked for an abelian subgroup $G\subset \SL_3(\C)$ in
\secref{sec:messy}.
Our proof below works in the
$2$-dimensional case with obvious modifications. And the vanishing
corresponding to \eqref{eq:mainassum} is already checked. Thus we have (1) and (2) also in $2$-dimensional case.

The statement~(1) was conjectured by Reid~\cite{Maizuru} based on the
corresponding result in the $2$-dimensional case proved by
Gonzales-Sprinberg and Verdier~\cite{GV}.  The
statement~(2) seems new even in dimension $2$.
We conjecture that the assumption \eqref{eq:mainassum} holds for any
finite subgroup $G$ of $\SL_3(\C)$. Note that our statement makes
sense in any dimension provided $X$ is nonsingular, and we conjecture
it holds under a reasonable, yet unknown, assumption on $G$. Remark
that in dimension $4$, when $G$ is the group of order $2$ generated by
$\operatorname{diag}(-1,-1,-1,-1)$, the statement of \thmref{thm:main}
is false while $X$ is nonsingular. In this example, $X$ is not
crepant, and the complex is not exact. Thus the smoothness of $X$ and the condition~\eqref{eq:mainassum} are not equivalent at least in dimension $4$.

\thmref{thm:main} has many interesting applications. First, we prove
$K_X = \shfO_X$ (see \thmref{thm:crepant}). When $G$ is an abelian
subgroup, this was proved by Nakamura~\cite{Iku} using the description
of $X$ as a toric variety. Our proof uses only the above mentioned
vanishing of certain homology groups.

The second application is much more interesting. We consider the
intersection pairing $(\ ,\ )$ on $K^c(X)$ defined by
\begin{equation}\label{eq:intpair}
  (S, T) = \langle \theta S, T\rangle,
\end{equation}
where $\theta\colon K^c(X) \to K(X)$ is the natural homomorphism.
Then, we have the following relation between the intersection pairing
on $K^c(X)$ and the decomposition of the tensor product.

\begin{Corollary}\label{cor:int}
Assume the same assumption as in \thmref{thm:main}.
The intersection pairing on $K^c(X)$ and the tensor product
decomposition~\eqref{eq:tensor} are related by
\begin{equation*}
  ( S_k^\vee, S_l) = a^{(2)}_{kl} - a^{(1)}_{kl} = a^{(1)}_{lk} - a^{(1)}_{kl},
\end{equation*}
where $S_k^\vee$ is the dual of $S_k$:
\begin{equation*}
S_k^\vee = - \left[\tautR_k \longrightarrow \bigoplus_l a^{(1)}_{kl}
  \tautR_l \longrightarrow \bigoplus_l a^{(2)}_{kl} \tautR_l
  \longrightarrow \tautR_k\right],
\end{equation*}
and the second equality follows from ${\textstyle\bigwedge^2} Q =
Q^*$.
\end{Corollary}

This corollary follows from $\theta S_k^\vee = \sum_l (a^{(2)}_{kl} -
a^{(1)}_{kl})\tautR_l$ and the above theorem.

In dimension $2$, the corresponding statement turns out to be
\begin{equation*}
(S_k^\vee, S_l) = 2\delta_{kl} - a^{(1)}_{kl}.
\end{equation*}
In \secref{sec:2d}, we will express $S_k$ in terms of irreducible
components $C_l$. In this way, we recover the identification of the
decomposition of tensor products and the intersection pairing
explained in the introduction.

Unfortunately we could not give an explicit expression of $S_k^\vee$
in the linear combination of $S_l$'s (or equivalently $\tautR_k^\vee$
in $\tautR_l$'s) in general. Thus we could not determine the
intersection product in terms of $G$.

Now we begin the proof of \thmref{thm:main}. First we show
\begin{Theorem}\label{thm:span}
$\{ \tautR_k \}_{k=0}^r$ and $\{ S_k \}_{k=0}^r$ span
$K(X)$ and $K^c(X)$ respectively.
\end{Theorem}

\begin{proof}
Modifying \eqref{eq:kosvec}, we introduce the following complex
depending on a parameter $s$:
\begin{equation}\label{eq:defkos}
\begin{CD}
  C_s\colon  F_3 @>d_3^s>> F_2 @>d_2^s>> F_1 @>d_1^s>> F_0,
\end{CD}
\end{equation}
where $F_k$ is as in \eqref{eq:kosvec} and $d_k^s$ is given by
\begin{equation*}
  d_k^s (\eta) \defeq s B^2\wedge \eta - \eta\wedge B^1\quad
  \text{for $\xi\in F_k$}.
\end{equation*}
This is still complex thanks to the equation $[B^a\wedge B^a] = 0$ ($a
=1, 2$). When $s = 1$, it is the original complex, which gives us a
resolution of $\shfO_{\Delta}$. If $s \neq 0$, this is nothing but the
pull-back of the complex $C_{s=1}$ by the automorphism of $X\times X$
defined by
\begin{equation}\label{eq:action}
  ([B^1, i^1], [B^2, i^2]) \mapsto ([B^1, i^1], [sB^2, i^2]).
\end{equation}

When $s = 0$, the complex~\eqref{eq:defkos} decomposes as
\begin{equation*}
\bigoplus_k p_1^*S_k \otimes p_2^* \tautR_k,
\end{equation*}
where $p_i$ is the projection to the $i$th factor.

Let $\Supp C_s$ be the subvariety on which the complex $C_s$ is not
exact. When $s = 1$, $\Supp C_s$ is the diagonal. For $s\neq 0$,
$\Supp C_s$ is the pull-back of the diagonal by \eqref{eq:action}. For
$s = 0$, $\Supp C_{s=0}$ is contained in $\pi^{-1}(0)\times X$ by
\lemref{lem:support}. In particular, in each case, the restriction of
the first projection $p_2\colon\Supp C_s \to X$ is proper. We consider
$C_s$ as a complex on $X\times X\times \C$, pulling back vector
bundles $F_k$ and setting the differential $d_k^s$ on $X\times X\times
\{ s\}$.  Then we can define the operator by
\begin{equation*}
  K(X) \ni E \longmapsto p_{23*}(p_1^* E \otimes C_s)\in K(X\times \C),
\end{equation*}
where $p_{23}\colon X\times X\times\C \to X\times\C$ is the projection 
to the second and the third factor.

Let $p\colon X\times\C \to X$ be the projection. It is known that
$p^*\colon K(X)\to K(X\times\C)$ is an isomorphism~\cite[IX,
1.6]{SGA6}.  Let $a_s\colon X\to X\times \C$ denote the embedding
given by $a_s(x) = (x,s)$. It satisfies $a_s^* p^* = (p\circ a_s)^* =
\operatorname{id}$, and hence $a_s^*$ is independent of $s$. If we
choose $s = 1$, we have $E = a_1^* p_{23*}(p_1^* E \otimes C_s)$ since
$C_{s=1}$ is the resolution of the diagonal. Comparing with the
pull-back by $a_0$, we get
\begin{equation}\label{eq:subst}
  E = a_0^* p_{23*}(p_1^* E\otimes C_s)
      = \sum_k \langle E, S_k\rangle \tautR_k ,
\end{equation}
where $\langle\ ,\ \rangle$ is the pairing given by
\eqref{eq:pairing}. In particular, this shows that $\{\tautR_k\}$
generates $K(X)$ as $\mathbb Z$-modules.

Similarly, we consider 
\begin{equation*}
  p_{13*}(p_2^* S\otimes C_s)
\end{equation*}
for $S\in K^c(X) = K(\pi^{-1}(0))$. Since $p_2\colon\Supp C_s \to X$
is proper, this defines an operator from $K(\pi^{-1}(0))$ to
$K(\pi^{-1}(0)\times\C)$. The pull-back homomorphism is
independent of $s$ as above, hence
\begin{equation*}
  S = a_0^*p_{13*}(p_2^* S\otimes C_s)
    = \sum_k \langle \tautR_k, S\rangle S_k.
\end{equation*}
This implies that $\{ S_k \}$ generates $K^c(X)$.
\end{proof}

We postpone the proof of the linear independence of $\tautR_k$ until
the end of this section.

The following is the first application of \thmref{thm:span}
\begin{Theorem}
$X$ is connected.
\end{Theorem}

\begin{proof}
Let $\{ X_\alpha \}$ be the set of connected components of $X$. For
any locally free sheaf $E$ (whose rank may change on components), we
assign the rank of its restriction to $X_\alpha$. Then we define the
{\it augmentation}
\begin{equation*}
  \varepsilon\colon K(X) \to {\mathbb Z}^{\pi_0(X)},
\end{equation*}
where $\pi_0(X)$ denotes the set of connected components of $X$. This
is surjective. However, $K(X)$ is generated by the tautological vector
bundles $\tautR_k$'s which have constant rank over the whole $X =
\coprod X_\alpha$. Hence $X$ must be connected.
\end{proof}

As we promised above, we prove that the canonical bundle $K_X$ is
trivial as an application of \thmref{thm:span}. Another ingredient is
the Serre duality. For an element $S\in K^c(X)$ represented by a
complex
\begin{equation*}
E_n \xrightarrow{d_n} E_{n-1} \xrightarrow{d_{n-1}} \cdots
\xrightarrow{d_2} E_1 \xrightarrow{d_1} E_0,
\end{equation*}
we define its dual $S^\vee\in K^c(X)$ by
\begin{equation*}
(-1)^n \left[ E_0^\vee \xrightarrow{{}^td_0} E_1^\vee
\xrightarrow{{}^td_1} \cdots \xrightarrow{{}^td_{n-1}} E_{n-1}^\vee
\xrightarrow{{}^td_n} E_n^\vee \right].
\end{equation*}

Then the Serre duality implies that
\begin{equation}\label{eq:Serre}
\langle E, S\rangle = - \langle E^\vee\otimes K_X, S^\vee\rangle
\end{equation}

\begin{Theorem}\label{thm:crepant}
The canonical bundle $K_X$ is trivial in $\operatorname{Pic}(X)$.
\end{Theorem}

\begin{proof}
Since the composition $\operatorname{Pic}(X)\to K(X)
\xrightarrow{\det}\operatorname{Pic}(X)$ is the identity, it is enough
to show that $K_X = \shfO_X$ in $K(X)$.

Substituting $E = \tautR_0$, $S = S_k^\vee$ into \eqref{eq:Serre}, we have
\begin{equation}\label{eq:dual}
 \langle \tautR_0,  S_k^\vee\rangle = - \langle K_X, S_k\rangle 
\end{equation}
where we have used $\tautR_0 = \tautR_0^\vee = \shfO_X$. Combining
with \eqref{eq:subst}, we get
\begin{equation}\label{eq:cano}
  K_X = - \sum_k \langle \tautR_0,  S_k^\vee\rangle \tautR_k.
\end{equation}

On the other hand, if we replace
$d^s_k$, in the proof of \thmref{thm:span}, by
\begin{equation*}
  d_k^{\prime s} (\eta) \defeq B^2\wedge \eta - s\eta\wedge B^1,
\end{equation*}
we get a homotopy between the complex~\eqref{eq:kosvec} and
\begin{equation*}
  \bigoplus_k p_1^*\tautR_k^\vee \otimes p_2^*(- S_k^\vee).
\end{equation*}
By the same argument as above, we obtain
\begin{equation*}
  E = -\sum_k \langle E, S_k^\vee\rangle \tautR_k^\vee
\end{equation*}
instead of \eqref{eq:subst}. Applying the duality, we have
\begin{equation*}
  E^\vee = -\sum_k \langle E, S_k^\vee\rangle \tautR_k.
\end{equation*}
Substituting $E = \tautR_0$, we have
\begin{equation*}
  \tautR_0 = -\sum_k \langle \tautR_0, S_k^\vee\rangle \tautR_k.
\end{equation*}
Comparing with \eqref{eq:cano}, we get
\(
  K_X = \tautR_0 = \shfO_X.
\)
\end{proof}

Thus we get
\begin{Corollary}
$X$ is a crepant resolution of $\C^3/G$.
\end{Corollary}

Now we relate the representation ring and the cohomology group. Let
$\ch\colon K(X)\to H^*(X,\Q)$ be the Chern character homomorphism. As
we used $K^c(X)$ besides $K(X)$, we need to consider the cohomology
group with compact support $H^*_c(X,\Q) \cong
H^*(X,X\setminus\pi^{-1}(0),\Q)$. The isomorphism can be shown by
observing the $\C^*$-action induced by $(x,y,z)\mapsto (tx, ty, tz)$
retracts $X$ to a neighborhood of $\pi^{-1}(0)$. We have the localized
Chern character homomorphism $\ch^c\colon K^c(X)\to H^*_c(X,\Q)$
defined by Iversen~\cite{Ive}.

\begin{Theorem}
\rom{(1)} The rational cohomology groups $H^*(X,\Q)$, $H^*_c(X,\Q)$
vanish in odd degrees.

\rom{(2)} $\{ \ch(\tautR_k) \}_{k=0}^r$ and $\{ \ch^{c}(S_k)
\}_{k=0}^r$ form dual bases of $H^*(X,\Q)$, $H^*_c(X,\Q)$ with respect
to the pairing
\begin{equation*}
\langle \alpha,\beta\rangle
= \int_X \alpha\cup\beta \operatorname{Td}(X),
\end{equation*}
where $\operatorname{Td}(X)$ is the Todd class of $X$.
\end{Theorem}

\begin{proof}
The proof proceed exactly as \thmref{thm:span} if we apply either the
usual Chern character or the localized Chern character. Then we find
that $\{ \ch(\tautR_k) \}_{k=0}^r$ and $\{\ch^{c}(S_k)\}_{k=0}^r$ span
$H^*(X,\Q)$, $H^*_c(X, \Q)$ respectively. Thus we have the
assertion~(1). Moreover, we know $\dim H^*(X,\Q) = r+1$ by previous
results on McKay correspondence \cite{BD,IR}. Hence they are
bases. Substituting $E = \tautR_l$ into \eqref{eq:subst}, applying the
Chern character, and using the Riemann-Roch, we get
\begin{equation*}
  \ch(\tautR_l) = \sum_k \langle \tautR_l, S_k\rangle \ch(\tautR_k)
  = \sum_k \langle \ch(\tautR_l), \ch^{c}(S_k)\rangle 
   \ch(\tautR_k).
\end{equation*}
This shows that $\langle \ch(\tautR_l), \ch^{c}(S_k)\rangle
= \delta_{kl}$. Thus we have the assertion~(2).
\end{proof}

In the course of the proof, we proved the linear independence of
$\tautR_k$ and the equality $\langle\tautR_l, S_k\rangle =
\delta_{kl}$. Thus we have completed the proof of
\thmref{thm:main}.

\section{$2$-dimensional case}\label{sec:2d}

In this section, we study $\tautR_k$ and $S_k$ in more detail in the
$2$-dimensional case. Since this case was already studied before
\cite{IN,Na-quiver}, we only give a sketch.

The same argument as in the previous section shows that $\{ \tautR_k
\}$ and $\{ S_k\}$ form dual base of $K(X)$ and $K^c(X)$, where $S_k$
in this case is
\begin{equation}\label{eq:2dSk}
S_k\colon \tautR_k^\vee \xrightarrow{{}^t(\wedge B)}
\bigoplus_l a_{kl} \tautR_l^\vee 
\xrightarrow{{}^tB} \tautR_k^\vee.
\end{equation}

First suppose $k\neq 0$. Then the cyclic vector condition implies that
${}^t(\wedge B)$ is injective on each fiber. The other homology
groups can be determined by using the following.

\begin{Proposition}\label{prop:2dcoh}
\rom{(1)} Let $C_k$ denote the subvariety where ${}^tB$ in
\eqref{eq:2dSk} is not surjective. Then $C_k$ is isomorphic to the
complex projective line $\CP^1$.

\rom{(2)} On $C_k$, the cokernel of ${}^tB$ is isomorphic to
$\shfO_{C_k}(-1)$. \rom(Here $-1$ means the dual of the hyperplane
bundle on $C_k\cong\CP^1$.\rom)

\rom{(3)} The restriction of the tautological bundle $\tautR_l$ to $C_k$ is
\begin{equation*}
\begin{cases}
  \shfO_{C_k}(1)\oplus (\shfO_{C_k})^{\oplus\rank\tautR_l -1}
    & \text{if $k = l$}, \\
  (\shfO_{C_k})^{\oplus\rank\tautR_l}
    & \text{if $k \neq l$}.
\end{cases}
\end{equation*}
\end{Proposition}

Thus the homology of $S_k$ vanishes in degree $1$, $2$ and
$\shfO_{C_k}(-1)$ in degree $0$. Hence $S_k \cong
\shfO_{C_k}(-1)$. The equality $\langle \tautR_l, S_k\rangle =
\delta_{kl}$ follows directly in this case.

Next suppose $k = 0$. Instead of considering $S_0$, we study $S_0^\vee$:
\begin{equation}\label{eq:2dS0}
S_0\colon \tautR_0 \xrightarrow{B}
\bigoplus_l a_{0l} \tautR_l
\xrightarrow{\wedge B} \tautR_0.
\end{equation}
Then $B$ is injective by the cyclic vector condition.
\begin{Proposition}\label{prop:2dcoh2}
\rom{(1)} The homomorphism $\wedge B$ in \eqref{eq:2dS0} is not
surjective exactly on the exceptional set $\pi^{-1}(0)$ of the
resolution $X\to \C^2/G$.

\rom{(2)} On the exceptional set, the cokernel of $\wedge B$ is
$\shfO_{\pi^{-1}(0)}$.
\end{Proposition}

Since $\langle \tautR_k, S_0\rangle = \langle \tautR_k^\vee,
S_0^\vee\rangle$ by the Serre duality, we could check $\langle
\tautR_k, S_0\rangle = \delta_{k0}$ also in this case.

Using the above, it becomes easy to determine $\ch^{c}(S_k)$. First
notice that $\operatorname{Td}(X) = 1$ in this case. Recall also that
$H^4_c(X,\mathbb Q) = \mathbb Q \Omega$ where $\Omega$ is the
canonical generator satisfying $\int_X \Omega = 1$. Then $\langle
\tautR_0, S_k\rangle = \delta_{k0}$ implies that the degree $4$ part
of $\ch^{c}(S_k) = \delta_{k0}\Omega$ because $\ch(\tautR_0)$ is the
canonical generator of $H^0(X,\mathbb Q)$.  By \propref{prop:2dcoh},
the degree $2$ part of $\ch^{c}(S_k)$ for $k\neq 0$ is the Poincar\'e
dual of $[C_k]$. Since $\ch^{c}(S_k)$ is a basis of $H^*_c(X,\mathbb
Q)$, this implies that $C_k$'s are irreducible components of the
exceptional set. (This can be directly checked by studying $C_k$'s.)
Then the equality $\langle \tautR_k, S_l\rangle = \delta_{kl}$ means
that $\{c_1(\tautR_k)\}_{k\neq 0}$ is the dual basis of $\{ [C_k ]\}$.
This is the second statement of a result of Gonzales-Sprinberg and
Verdier explained in the introduction. Finally, considering
\begin{equation*}
\begin{split}
  0 &= \langle \tautR_k, S_0\rangle \qquad\text{for $k\neq 0$} \\
    &= \langle c_1(\tautR_k), \ch^{c}(S_0)\rangle + \rank \tautR_k,
\end{split}
\end{equation*}
we get
\begin{equation*}
  \text{degree $2$ part of }\ch^{c}(S_0)
   = -\sum_{k=1}^r \rank\tautR_k \operatorname{P.D.}[C_k],
\end{equation*}
where $\operatorname{P.D.}$ is the Poincar\'e dual.

Note also that 
\begin{equation*}
\ch^{c}(S_k^\vee) =
\begin{cases}
  - \ch^{c}(S_k) & \text{for $k\neq 0$}, \\
  \sum_{k=1}^r \rank\tautR_k \operatorname{P.D.}[C_k] + \Omega
   & \text{for $k = 0$}.
\end{cases}
\end{equation*}
This together with $(S_k^\vee, S_l) = (\ch^{c}(S_k^\vee), \ch^{c}(S_l))
= 2\delta_{kl} - a_{kl}$ determine the intersection pairing. In
particular, we have $(S_k, S_l) = (\operatorname{P.D.}[C_k],
\operatorname{P.D.}[C_l]) = a_{kl} - 2\delta_{kl}$ for $k,l\neq
0$. Thus we have checked the identification of the intersection pairing
and the decomposition of tensor products without using the
classification of simple singularities.

There are possibly many ways to prove
Propositions~\ref{prop:2dcoh},\ref{prop:2dcoh2}, since we have more
or less explicit description of $X$ (e.g.~\cite{IN}). However, we
would like to remark that these follows from the theory of the quiver
varieties introduced by the second author~\cite{Na-Duke} without the 
knowledge of the explict description.
The variety $X$ is an example of quiver varieties associated with the
extended Dynkin diagram. The second author defined the Hecke
correspondence in the product of quiver varieties [ibid., 10.4]. The
locus $C_k$ ($k\neq 0$) is an example of Hecke correspondence, where
the one factor is $X$ and the other factor is the scheme parametrising
$0$-dimensional subschemes $Z$ of $\C^2$ such that
\begin{enumerate}
\item $Z$ is invariant under the $G$-action,
\item $H^0(\shfO_Z) \oplus \rho_k$ is the regular representation of
$G$,
\end{enumerate}
which consists of one point. The assertions in \propref{prop:2dcoh}
follow from [ibid., Lemma~10.10 and its proof].
(We omit details.)
\propref{prop:2dcoh2} is, in fact, much easier to prove, and holds
also in $3$-dimensional case. Since $\tautR_0$ is the trivial rank $1$
bundle over $X$, $\wedge B$ is not surjective only when $\wedge B =
0$. If $\wedge B = 0$ at $Z\in X$, $B$ is nilpotent as an endomorphism
of $H^0(\shfO_Z) = A/I$, hence contained in $\pi^{-1}(0)$. Conversely
if $Z\in\pi^{-1}(0)$, we have a filtration on $H^0(\shfO_Z)$, under
which $B$ is strictly upper triangular by the Hilbert criterion as in
\cite{Na-affine}. It implies $\wedge B = 0$. 

\section{Toric resolution: the case $G\subset \SL_3(\C)$ abelian}
\label{sec:messy}

In this section, we assume $G$ is an abelian subgroup of
$\SL_3(\C)$. It was proved by Nakamura~\cite{Iku} that $X$ is a
crepant resolution of $\C^3/G$ under this assumption. (There is also
an explanation of Nakamura's proof by Reid~\cite{Maizuru}.)

\begin{Theorem}[Nakamura\cite{Iku}]\label{thm:toric}
If $G$ is a finite abelian subgroup of $\SL_3(\C)$, then $X$ is a
crepant resolution of $\C^3/G$.
\end{Theorem}

In this section, we prove the following which is stronger than the
smoothness of $X$:
\begin{Theorem}\label{prop:assum}
Under the same assumption as \thmref{thm:toric}, the following holds:
\begin{equation}\label{eq:assum}
 \dim \Hom_A(I_1,A/I_2)^G =
 \begin{cases}
 3, & \text{when $I_1 = I_2$,}\\
 1, & \text{when $I_1 \neq I_2$}.
 \end{cases}
\end{equation}
\end{Theorem}

The rest of this section is devoted to the proof of
\thmref{prop:assum}. Though the smoothness of $X$, i.e.\ 
\eqref{eq:assum} for $I_1 = I_2$, is contained in Nakamura's
result~\cite{Iku}, we give its proof for the sake of the reader. Note
also that the crepantness of $X$ follows from \thmref{thm:crepant} and 
the result in this section.
Our proof for the smoothness is almost the same as Nakamura's, and the
technique (e.g., the use of the diagram $\frak J$, \lemref{lem:ae}) is
due to him.

By changing coordinates, we may assume $G$ is diagonal.  We consider
an action of the three dimensional torus on $\C^3$ defined by $t \cdot
(x,y,z) =(t_1 x, t_2 y, t_3 z)$ for $t = (t_1, t_2, t_3)$.  It induces
an action on the Hilbert scheme, commuting with the action of
$G$. Thus this action induces the torus action on $X$.

In the sequel, the fixed points of the torus action will play the
crucial role. These correspond to ideals generated by monomials.
First we will classify all these ideals.  Then we will check
\eqref{eq:assum} for all these ideals. This implies \eqref{eq:assum}
for general $I_1$, $I_2$ as explained later.

\subsection{Classification of fixed points}
The inclusion $G \subset \SL_3(\C)$ induces an action of $G$ on the
coordinate ring $A=\C[ x,y,z]$. It decomposes into the sum of
irreducible representations, which are of the forms $\C x^ly^mz^n$ for
some $l$, $m$, $n\geq 0$. If $I \in X$, then $A/I$ is the regular
representation. Hence each irreducible representation appears in $A/I$
with multiplicity one. This implies the following very useful lemma,
which will be used throughout in this section.

\begin{Lemma}\label{lem:rep}
Let $I \in X$. Suppose that $x^ly^mz^n$ and $x^{l'}y^{m'}z^{n'}$ are
two different monomials which give the isomorphic irreducible
representation of $G$. Then at least one of them is contained in $I$.
\end{Lemma}
 
For the study of ideals, we use the following graphical description of
monomials in $A/(xyz)$:
\begin{equation*}
\renewcommand{\arraystretch}{1.4}
\renewcommand{\arraycolsep}{8pt}
\setcounter{MaxMatrixCols}{13}
\begin{matrix}
 & & & & & & & & & & & &\\[-18pt]
 & \kern-1cm \cdots \kern-1cm \\
 && \kern-1cm y^2 \kern-1cm \\
 &&& \kern-1cm y \kern-1cm && \kern-1cm xy \kern-1cm &&
 \kern-1cm x^2y \kern-1cm \\
\kern-1cm y^2z^2 \kern-1cm &&
\kern-1cm yz \kern-1cm &&
\kern-1cm 1 \kern-1cm &&
\kern-1cm x \kern-1cm &&
\kern-1cm x^2 \kern-1cm &&
\kern-1cm x^3 \kern-1cm &&
\kern-1cm \cdots \kern-1cm \\
 &&& \kern-1cm z \kern-1cm && \kern-1cm xz \kern-1cm \\
 && \kern-1cm \cdots \kern-1cm
\end{matrix}
\end{equation*}
Here $xyz$ is not drawn since it is always in ideals as it gives an
isomorphic representation as $1$.

\begin{Proposition}[cf.~\protect{\cite[7.2]{Maizuru}}]\label{prop:AB}
The ideal $I$ which is fixed by the torus action is written
as one of the following:
\begin{align}
\langle x^{a+d-1},\ y^{b+e-1},\ z^{c+f-1},\ x^ay^e,\ y^bz^f,\ z^cx^d,\ 
xyz \rangle \tag*{(A)}\label{eq:A}\\
\langle x^{a+d},\ y^{b+e},\ z^{c+f},\ x^ay^e,\ y^bz^f,\ z^cx^d,\ xyz \rangle
\tag*{(B)}\label{eq:B}
\end{align}
where $a,b,c,d,e,f > 0$.

Moreover, $x^{a+d-1}$ (resp. $x^{a+d}$) and $y^{b-1} z^{f-1}$ give the
isomorphic representation in type (A) (resp. (B)). Similar conditions
hold if we exchange $x$, $y$ and $z$.
\end{Proposition}

In above description, there are degenerate cases, for example, $e=1$
in (A) where $y^bz^f$ is not a generator.  But we could determine $a$,
$b$, $c$, $d$, $e$ and $f$ so that the conditions of the isomorphic
representations above hold (see \lemref{lem:ae}). We use the above
description even in degenerate cases as convention.

The proof of this proposition occupies this subsection. Let $I$ be an
ideal which is fixed by the torus action and generated by
monomials. By \lemref{lem:rep}, each irreducible representation of $G$
corresponds to the unique monomial in $I^c$.

Let $\alpha$, $\beta$, $\gamma$ be the exponents of the generators of 
$x^\bullet$, $y^\bullet$, $z^\bullet$
respectively. 

\begin{Lemma}\label{lem:generator}
Suppose $x^ay^e\in I$, $x^{a-1}y^e\notin I$, and $x^ay^{e-1}\notin I$
for $a, e \geq 1$, i.e.\  $x^ay^e$ is
a generator of $I$.  Then the unique monomial in $I^c$ which
has the isomorphic representation as $x^ay^e$ is $z^{\gamma -1}$.
\end{Lemma}
\begin{proof}
Let us consider the monomial $\varphi$ in $I^c$ which gives the isomorphic
representation of $G$ as $x^ay^e$.

If $\varphi =x^py^q$ with $p\geq 1$, then $x^{p-1}y^q$ gives the isomorphic
representation as $x^{a-1}y^e$. From the assumption $x^{a-1}y^e
\not\in I$ and \lemref{lem:rep}, we have $x^{p-1}y^q \in
I$. Therefore $x^py^q \in I$. This is a contradiction.  Exchanging $y$
and $z$, we can also eliminate the case $\varphi =x^pz^r$ for $p \geq 1$.

Next assume $\varphi =y^qz^r$ with $q \geq 1$. Then $y^{q-1}z^r$ gives the
isomorphic representation as $x^ay^{e-1}$. Then we have $y^{q-1}z^r
\in I$ by Lemma~\ref{lem:rep}. This means that $y^qz^r \in I$ and
contradicts with the definition of $\varphi$. Thus we have $\varphi=z^r$. Since
$z^{r+1}$ gives the isomorphic representation as $x^ay^e$, we have
$z^{r+1} \in I$ by \lemref{lem:rep}.  Hence $z^{r +1}$ must be a
generator of $I$, i.e.\  $z^r = z^{\gamma -1}$.
\end{proof}

This lemma implies that $I$ has at most one generator of the form
$x^ay^e$ $(a,e\geq 1)$: If $x^ay^e$ and $x^{a'}y^{e'}$ are generators,
both $x^{a-1}y^{e-1}$ and $x^{a'-1}y^{e'-1}$ give the isomorphic
representation as the generator $z^{\gamma}$. This contradicts with
\lemref{lem:rep}.

Although we assume $x^ay^e$ is a generator in the lemma, we can show
the following even if there is no generator of that form (i.e.\ 
degenerate case).

\begin{Lemma}\label{lem:ae}
The monomial in $I^c$ which has the isomorphic representation as $z^{\gamma}$ 
is of the form $x^{a-1}y^{e-1}$ for some $a,e\geq 1$.
\end{Lemma}

This lemma can be proved as \lemref{lem:generator}. So we omit the proof.

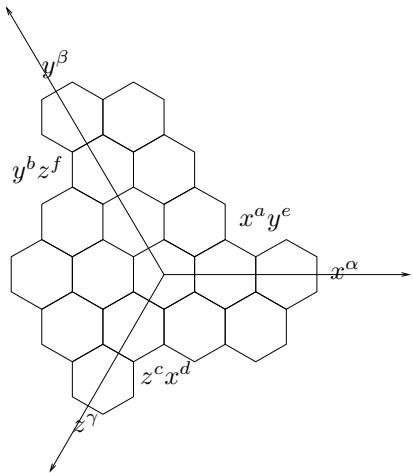
\begin{figure}[htbp]
\setlength{\unitlength}{0.0003in}
\begin{center}
\begin{picture}(6924,7989)(0,-10)
\path(2187,3087)(2187,3687)(1662,3987)
        (1137,3687)(1137,3087)(1662,2787)(2187,3087)
\path(2187,1287)(2187,1887)(1662,2187)
        (1137,1887)(1137,1287)(1662,987)(2187,1287)
\path(5337,3087)(5337,3687)(4812,3987)
        (4287,3687)(4287,3087)(4812,2787)(5337,3087)
\path(4812,2187)(4812,2787)(4287,3087)
        (3762,2787)(3762,2187)(4287,1887)(4812,2187)
\path(3237,3087)(3237,3687)(2712,3987)
        (2187,3687)(2187,3087)(2712,2787)(3237,3087)
\path(4287,3087)(4287,3687)(3762,3987)
        (3237,3687)(3237,3087)(3762,2787)(4287,3087)
\path(2712,2187)(2712,2787)(2187,3087)
        (1662,2787)(1662,2187)(2187,1887)(2712,2187)
\path(1662,5787)(1662,6387)(1137,6687)
        (612,6387)(612,5787)(1137,5487)(1662,5787)
\path(2712,5787)(2712,6387)(2187,6687)
        (1662,6387)(1662,5787)(2187,5487)(2712,5787)
\path(3762,3987)(3762,4587)(3237,4887)
        (2712,4587)(2712,3987)(3237,3687)(3762,3987)
\path(3237,4887)(3237,5487)(2712,5787)
        (2187,5487)(2187,4887)(2712,4587)(3237,4887)
\path(1137,3087)(1137,3687)(612,3987)
        (87,3687)(87,3087)(612,2787)(1137,3087)
\path(1662,3987)(1662,4587)(1137,4887)
        (612,4587)(612,3987)(1137,3687)(1662,3987)
\path(1662,2187)(1662,2787)(1137,3087)
        (612,2787)(612,2187)(1137,1887)(1662,2187)
\path(2187,4887)(2187,5487)(1662,5787)
        (1137,5487)(1137,4887)(1662,4587)(2187,4887)
\path(3762,2187)(3762,2787)(3237,3087)
        (2712,2787)(2712,2187)(3237,1887)(3762,2187)
\path(2712,3387)(6912,3387)
\path(6792.000,3357.000)(6912.000,3387.000)(6792.000,3417.000)
\path(2712,3387)(12,7962)
\path(98.827,7873.903)(12.000,7962.000)(47.154,7843.407)
\path(2712,3387)(762,12)
\path(796.057,130.912)(762.000,12.000)(848.009,100.895)
\put(100,4950){\makebox(0,0)[lb]{$y^bz^f$}}
\put(612,6762){\makebox(0,0)[lb]{$y^\beta$}}
\put(4000,4137){\makebox(0,0)[lb]{$x^ay^e$}}
\put(5562,3312){\makebox(0,0)[lb]{$x^\alpha$}}
\put(2300,1512){\makebox(0,0)[lb]{$z^cx^d$}}
\put(1137,687){\makebox(0,0)[lb]{$z^\gamma$}}
\end{picture}
\caption{Junior diagram}
\label{fig:junior}
\end{center}
\end{figure}

By above lemmas, the ideal is generated by $x^\alpha$, $y^\beta$,
$z^\gamma$, $x^ay^e$, $y^bz^f$, $z^cx^d$ and $xyz$ where $\alpha$,
$\beta$, $\gamma$, $a$, $b$, $c$, $d$, $e$, $f > 0$. We draw a diagram
by placing hexagons at the monomials in $I^c$.  (See
figure~\ref{fig:junior}.)  Let us call this {\it junior diagram}
$\frak J $. (This is called a $G$-graph in \cite{Iku}.)
We may also write the representations corresponding to the
monomials. In dimension 2, we had a similar correspondence between
idelas and Young diagrams (cf.\ \cite[Chapter 7]{Lecture}).

\begin{Definition}
Let $\varphi$ be a monomial in $A/(xyz)$. Let $\varphi_0$ be the
monomial in $I^c$ which has isomorphic representation of $G$ as
$\varphi$. (It exists and is unique since $A/I$ is the regular
representation.) We move the junior diagram $\frak J$ by the parallel
transport which maps $\varphi_0$ to $\varphi$, and denote the
transported diagram by ${\frak J}(\psi)$. We also call ${\frak
J}(\varphi)$ the {\it junior diagram}.
\end{Definition}

The parallel transport respects representations corresponding to
monomials. Hence each representation appears at the same position in
each junior diagram. Each monomials in $\varphi\in A/(xyz)$ belongs to
the unique junior diagram ${\frak J}(\varphi)$. Thus each junior
diagram is a kind of `fundamental domain' of $A/(xyz)$.

To complete the proof of \propref{prop:AB}, we have to show the
relation between the exponents of the generators. First we will see
the following:

\begin{Lemma}
We have one of the following:
\begin{subequations}
\begin{alignat}{3}
& \alpha\le a+d-1,&\quad 
     &\beta \le b+e-1,&\quad
     &\gamma \le c+f-1,\label{eq:one} \\
& \alpha \ge a+d, &\quad
     &\beta \ge b+e,& \quad
     &\gamma \ge c+f.\label{eq:two}
\end{alignat}
\end{subequations}
\end{Lemma}
\begin{proof}
{\bf Step 1}. Assume
\begin{subequations}
\begin{align}
& b+e \leq \beta, \label{eq:eins}\\
\intertext{and}
& \alpha \leq a+d-1. \label{eq:zwei} 
\end{align}
\end{subequations}
Since $x^a y^e$ has the isomorphic representation as $z^{\gamma-1}$,
$x^\alpha y^e$ has the isomorphic representation as
$z^{\gamma-1}x^{\alpha-a}$, which is in ${\frak J} = {\frak J}(1)$ by
\eqref{eq:zwei}.

On the other hand, $x^\alpha$ has the isomorphic representation as
$y^{b-1} z^{f-1}$. Hence $x^\alpha y^e$ has the isomorphic
representation as $y^{b+e-1} z^{f-1}$, which is in ${\frak J}(1)$ by
\eqref{eq:eins}. As $b,e\ge 1$, $y^{b+e-1} z^{f -1} \neq x^{\alpha -a}
z^{\gamma -1}$. Therefore there are two monomials in $\frak J (1)$
which has the isomorphic respresentations as $x^\alpha y^e$. This
contradicts with \lemref{lem:rep}. Hence we have either
\begin{subequations}
\begin{align}
& b+e -1 \ge \beta, \label{eq:n-eins}\\
\intertext{or}
& \alpha \ge a+d. \label{eq:n-zwei} 
\end{align}
\end{subequations}

{\bf Step 2}. Suppose $b+e-1\ge \beta$. Exchanging $x$ and $y$ in
step~1, we have either
\begin{align*}
& a+d -1 \ge \alpha, \\
\intertext{or}
& \beta \ge b+e. 
\end{align*}
By the assumption, the second case does not occur. Thus we have
$a+d-1\ge \alpha$. Exchanging $x$ and $z$, we also have
$c+f-1\ge\gamma$, and hence \eqref{eq:one}. Similar argument shows
that $\alpha\ge a+d$ implies \eqref{eq:two}.
\end{proof}

Now we complete the proof of \propref{prop:AB}. First we consider the
case \eqref{eq:one} and let us show that
\begin{equation}\label{eq:ichi}
 x^ay^\beta = x^a y^{b+e-1}.
\end{equation}
Since we have $\beta \leq b+e-1$, it is enough to show that the condition
\begin{equation}\label{eq:san}
\beta < b+e-1
\end{equation}
leads to a contradiction. 

Let us study the monomials near $x^{a-1}y^\beta$. 
Consider the parallel transport mapping $x^{a-1}y^\beta$ to
$y^{\beta-e+1}z^{\gamma}$, which respects representations of $G$. It maps
$x^a y^{\beta+1}$ to $ y^{\beta-e+1}z^{\gamma -1}$ which is in 
${\frak J}(1)$ by the assumption~\eqref{eq:san}. The monomial
$x^a y^\beta$ is also mapped to a monomial $y^{\beta -e} z^{\gamma -1}\in {\frak J}(1)$. 

On the other hand, $y^{\beta-e+1}z^{\gamma}$ is not in ${\frak J}(1)$.
Hence $x^a y^\beta$ and $x^a y^{\beta+1}$ are not in
$\frak J (x^{a-1}y^\beta)$. Since $x^{a-1}y^{\beta-1}\in{\frak J}(1)$, it
is not in $\frak J (x^{a-1}y^\beta)$ either.
In summary, for $x^{a-1}y^\beta$, the monomials in the upper right,
the right and the lower right are not in ${\frak
J}(x^{a-1}y^\beta)$. If $\varphi$ is the monomial in ${\frak J}(1)$, which
has the isomorphic representation as $x^{a-1} y^\beta$, it also has
the same property. Hence $\varphi$ must be $x^{\alpha -1}$.

Thus, $x^\alpha$ has the isomorphic representation as $x^a y^\beta$,
and hence as $y^{\beta -e}z^{\gamma-1}$, which is in ${\frak J}(1)$
by \eqref{eq:san}. However, $x^\alpha$ has the isomorphic
representation as $y^{b-1}z^{f-1}\in{\frak J}(1)$, hence we must have
$\gamma = f$, $\beta-e = b-1$ by \lemref{lem:rep}. But the latter
contradicts with the assumption~\eqref{eq:san}.

Therefore the condition \eqref{eq:ichi} holds in this
case. Exachanging $x$, $y$ and $z$, we have
\begin{equation}\label{eq:typeA}
\alpha = a+d-1, \quad \beta = b+e-1, \quad \gamma = c+f-1,
\end{equation}
i.e.\  the ideal is of type (A).

Next consider the case \eqref{eq:two}.
We will show that
\begin{equation}\label{eq:ni}
x^ay^\beta = x^{\alpha -d} y^\beta.
\end{equation}
Since we already have $a \leq \alpha -d$, it is enough to show that
\begin{equation}\label{eq:yon}
a < \alpha -d 
\end{equation}
leads to a contradiction. We study the monomials near $x^ay^{\beta -
1}$.  Consider the parallel transport mapping $x^a y^{\beta-1}$ to
$x^{a+d}z^{c}\notin{\frak J}(1)$. It maps $x^{a+1}y^\beta$ to
$x^{a+d}z^{c-1}$, which is in ${\frak J}(1)$ by the
assumption~\eqref{eq:yon}. Thus $x^{a+1}y^\beta$ (and hence $x^a
y^\beta$ also) is not in ${\frak J}(x^a y^{\beta-1})$. We have $x^a
y^{\beta-1}\notin{\frak J}(1)$ thanks to $\beta - 1 \ge b+e-1 \ge e$,
while we have $x^{a-1}y^{\beta-1}\in{\frak J}(1)$. Hence $x^{a-1}
y^{\beta-1}$ is not in ${\frak J}(x^a y^{\beta-1})$ either. By the
same argument as above, we conclude that the monomial in ${\frak
J}(1)$ which has isomorphic representation as $x^a y^{\beta-1}$ must
be $y^{\beta-1}$. Therefore $x^{a-1} y^{\beta-1}$ has the isomorphic
representation as $y^\beta z$, and hence as $z^c x^{d-1}$.

However, as $c < \gamma$ by \eqref{eq:two}, we have $z^c
x^{d-1}\in{\frak J}(1)$. Thus both $x^{a-1}y^{\beta -1}$ and
$z^cx^{d-1}$ are in $\frak J (1)$, and
$$ x^{a-1} y^{\beta -1} \neq z^c x^{d-1}$$
as $c\geq 1$. 
This is a contradiction. Thus we have \eqref{eq:ni}.
Exchanging $x$, $y$ and $z$, we have 
\begin{equation}\label{eq:typeB}
\alpha = a+d, \quad \beta = b+e, \quad \gamma = c+f,
\end{equation}
i.e.\  it is of type (B).

Combining the above two cases, we have only two types (A) and (B) for the generator of the ideal $I$ of a fixed point in $X$. 
This completes the proof of \propref{prop:AB}.

\subsection{Smoothness}
Using the above description, we obtain the smoothness of $X$:
\begin{Theorem}\label{thm:smooth}
$X$ is nonsingular of dimension $3$.
Moreover, $X$ is irreducible.
\end{Theorem}

This subsection is devoted to the proof of this theorem.

\begin{Lemma}
At every fixed point $I$, the dimension of the Zariski tangent space
is three, i.e.\ 
$$\dim\Hom_A(I, A/I)^G =3.$$
\end{Lemma}

\begin{proof}
Let $\Phi$ be a $G$-equivariant $A$-homomorphism $\Phi\colon I
\rightarrow A/I$. Let $\alpha$, $\beta$, $\gamma$ be the exponents of
the generators of $x^\bullet$, $y^\bullet$, $z^\bullet$ as before.

Let us consider the image of generators of $I$:
\begin{alignat*}{3}
& \Phi (x^\alpha)= p y^{b-1}z^{f-1} \bmod I,& \quad 
& \Phi (x^ay^e)=s z^{\gamma-1} \bmod I, &&\\
& \Phi (y^\beta)= q z^{c-1}x^{d-1} \bmod I,& \quad 
& \Phi (y^bz^f)=t x^{\alpha-1} \bmod  I,&\quad
& \Phi(xyz)=v \bmod I,\\
& \Phi (z^\gamma)= r x^{a-1}y^{e-1} \bmod I,& \quad 
& \Phi (z^cx^d)=u y^{\beta-1} \bmod I,&&
\end{alignat*}
where $p,q,r,s,t,u,v\in\C$.
Here we determine the image so that it has the isomorphic
representation of $G$ as the generator.

First suppose $I$ is of type (A), i.e.\  $\alpha = a + d -1$, $\beta =
b+e-1$, $\gamma = c+f-1$. Let us consider the image of 
$x^\alpha y^e = x^{a+d-1}y^e$.
We have 
$$
\Phi(x^{a+d-1}y^e)= p y^{b+e-1}z^{f-1}\bmod I= 0
$$
as $y^{b+e-1}z^{f-1}\in I$. On the other hand, we have
$$
\Phi(x^{a+d-1}y^e)=sx^{d-1}z^{c+f -2}\bmod I.
$$
Since $x^{d-1}z^{c+f -2}\not\in I$, we get $s=0$. Exchanging $x$,
$y$ and $z$, we obtain $t=u=0$. Considering the image of $x^a y^e z$
in two ways, we similarly get $v=0$. Therefore the dimension of the
Zariski tangent space at $I$ is three.

Next suppose $I$ is of type (B), i.e.\  $\alpha = a+d$, $\beta=b+e$,
$\gamma = c+f$. Then we have
$$
\Phi(x^{a+d}y^e)= p y^{b+e-1}z^{f-1}\bmod I,
$$
and $y^{b+e-1}z^{f-1}\not\in I$.
But
$$\Phi(x^{a+d}y^e)=sx^{d}z^{c+f-1}\bmod I = 0$$
as $x^dz^{c+f -1}\in I$. 
Then we have $p=0$. After exchanging $x$, $y$ and $z$, $q=r=0$ holds.
We also get $v = 0$ as in the case (A). Thus the assertion holds in
any cases.
\end{proof}

\begin{Lemma}\label{lem:nsfixed}
$X$ is nonsingular at fixed points of the torus action.
\end{Lemma}

\begin{proof}
Let $I\in X$ be a fixed point of the torus action. First suppose $I$
is of type (A).
Let us consider the following defining equation:
\begin{alignat}{3}
& x^{a+d-1} = \lambda y^{b-1}z^{f-1},& \quad 
& x^ay^e=\lambda\mu z^{c+f -2},& &\notag\\
& y^{b+e-1} = \mu z^{c-1}x^{d-1},& \quad 
& y^bz^f=\mu\nu x^{a+d -2},& \quad &xyz=\lambda\mu\nu, 
\label{eq:chart}\\
& z^{c+f-1} = \nu x^{a-1} y^{e-1},& \quad 
& z^cx^d=\nu\lambda y^{b+e -2}.& &\notag
\end{alignat}
This equations determines an ideal in the neighborhood of the fixed
point. ($\lambda=\mu=\nu=0$ is the fixed point.) If $\lambda\mu\nu\neq
0$, we have
\begin{equation*}
  \# G = 4-2( a + b + c + d + e + f) + ab+ac+ae+bc+bf+dc+de+df+ef 
\end{equation*}
distinct solutions of the above equation. Hence it corresponds to a
$G$-orbit consisting of distinct $\# G$-points, and has the Zariski
tangent space of dimension $3$. In particular, $X$ is of dimension $3$
in the neighbourhood of $I$. Since the Zariski tangent space at $I$ is
of dimension $3$, it implies that $X$ is nonsingular at $I$.

Next suppose $I$ is of the type (B). We consider the defining equations
\begin{alignat}{3}
& x^{a+d} = \nu\lambda y^{b-1}z^{f-1},& \quad 
& x^ay^e=\lambda z^{c+f -1},& &\notag\\
& y^{b+e} = \lambda\mu z^{c-1}x^{d-1},& \quad 
& y^bz^f=\mu x^{a+d -1},& \quad &xyz=\lambda\mu\nu,\label{eq:chart2} \\
& z^{c+f} = \mu\nu x^{a-1} y^{e-1},& \quad 
& z^cx^d=\nu y^{b+e -1}.& &\notag
\end{alignat}
This equation has
\begin{equation*}
  \# G =1-( a + b + c + d + e + f) +  ab+ac+ae+bc+bf+dc+de+df+ef 
\end{equation*}
distinct solutions if $stu\neq 0$. The above argument shows that $X$
is nonsingular at $I$ in this case.
\end{proof}

\begin{proof}[Proof of \thmref{thm:smooth}]
We take a generic one-parameter subgroup $\lambda\colon \C^*\to
(\C^*)^3$ such that $\lambda(t)\to 0$ as $t\to 0$. For any $I\in X$,
$\lambda(t)^*I$ converges to a fixed point of the torus action. Thus
$X$ is nonsingular at $I$ by \lemref{lem:nsfixed}. Hence $X$ is
nonsingular everywhere. The argument also shows that each connected
component of $X$ contains fixed points. However, the fixed points are
contained the component containing $G$-orbits of distinct
points. Therefore, $X$ must be connected.
\end{proof}

Since $X$ is nonsingular and has an action of $3$-dimensional torus
with an open dense orbit, we have

\begin{Corollary}
$X$ is a toric variety.
\end{Corollary}

The coordinate neighbourhoods~\eqref{eq:chart},\eqref{eq:chart2} are
affine charts around fixed points. 

The fan corresponding to $X$ are described in \cite{Iku}. So we do not 
reproduce it here.

\subsection{The case when $I_1$ and $I_2$ are in a common affine chart}

Our remaining task is to check \eqref{eq:assum} for $I_1$ and $I_2$ are 
contained in different affine charts of \eqref{eq:chart} or \eqref{eq:chart2}.
In this subsection, we check it when $I_1$ and $I_2$ are in a
common affine chart given by \eqref{eq:chart}.  One can check in the
case when both are in a chart given by \eqref{eq:chart2} in a similar
way, so we omit the proof.

Suppose $I_1$ is given by \eqref{eq:chart} and $I_2$ is given also by
\eqref{eq:chart} with $\lambda$, $\mu$, $\nu$ are replaced by
$\lambda'$, $\mu'$, $\nu'$. ($a,b,c,d,e,f$ are common.) By the
assumption $I_1\neq I_2$, we have $(\lambda,\mu,\nu)\neq
(\lambda',\mu',\nu')$. We may assume $\lambda\neq \lambda'$.  Let
$\Phi\colon I_1\to A/I_2$ be a $G$-equivariant $A$-homomorphism. We
determine images of generators:
\begin{align*}
& \Phi (x^{a+d-1}-\lambda y^{b-1}z^{f-1})=
        p y^{b-1}z^{f-1} \bmod I_2,\\
& \Phi (y^{b+e-1}-\mu z^{c-1}x^{d-1})=
        q z^{c-1}x^{d-1} \bmod I_2,\\
& \Phi (z^{c+f-1}-\nu x^{a-1}y^{e-1})=
        r x^{a-1}y^{e-1} \bmod I_2,\\
& \Phi (x^ay^e-\lambda\mu z^{c+f-2})=s z^{c+f-2} \bmod I_2,\\
& \Phi (y^bz^f-\mu\nu x^{a+d-2})=t x^{a+d-2} \bmod  I_2,\\
& \Phi (z^cx^d-\nu\lambda y^{b+e-2})=u y^{b+e-2} \bmod I_2,\\
& \Phi(xyz-\lambda\mu\nu)=v \bmod I_2.
\end{align*}
Consider the image of
\begin{equation*}
 (z^cx^d-\nu\lambda y^{b+e-2})x^{a-1}
= (x^{a+d-1}-\lambda y^{b-1}z^{f-1})z^c 
+ \lambda(z^{c+f-1}-\nu x^{a-1}y^{e-1})y^{b-1}.
\end{equation*}
We have
\begin{align*}
& u x^{a-1}y^{b+e-2}\bmod I_2 =
p y^{b-1}z^{c+f-1}\bmod I_2 + \lambda r x^{a-1}y^{b+e-2}\bmod I_2\\
=\;& (\nu' p + \lambda r)x^{a-1}y^{b+e-2}\bmod I_2.
\end{align*}
Since $x^{a-1}y^{b+e-2}\notin I_2$, we have
\(
u = \nu' p + \lambda r.
\)
Exchanging $x$ and $z$, we get
\(
u = \lambda' r + \nu p.
\)
Thus we have
\(
r = (\nu - \nu')p/(\lambda - \lambda').
\)
Then exchanging $x$, $y$ and $z$, we get
\(
s = \mu' p + \lambda q,
\)
\(
t = \nu' q + \mu r,
\)
\(
q = (\mu-\mu')p/(\lambda - \lambda').
\)

Next consider the image of
\begin{align*}
 & (xyz-\lambda\mu\nu)x^{a+d-2}y^{e-1} \\
=\;& (x^{a+d-1}-\lambda y^{b-1}z^{f-1})y^{e}z \\
&\qquad + \lambda(y^{b+e-1}-\mu z^{c-1}x^{d-1})z^{f}
+ \lambda\mu(z^{c+f-1}-\nu x^{a-1}y^{e-1})x^{d-1}. 
\end{align*}
We get
\begin{align*}
& v x^{a+d-2}y^{e-1}\bmod I_2 \\
=\; & p y^{b+e-1}z^{f}\bmod I_2
  + \lambda q z^{c+f-1}x^{d-1}\bmod I_2
  + \lambda\mu r x^{a+d-2}y^{e-1}\bmod I_2\\
=\; & (\mu'\nu' p + \lambda \nu' q + \lambda\mu r)
x^{a+d-2}y^{e-1}\bmod I_2.
\end{align*}
Since $x^{a+d-2}y^{e-1}\notin I_2$, we have
\(
v = \mu'\nu' p + \lambda \nu' q + \lambda\mu r.
\)
Thus $q$, $r$, $s$, $t$, $u$, $v$ are determined by $p$. Hence
we have $\dim \Hom_A(I_1, A/I_2)^G = 1$ in this case.

\subsection{Reduction to the case when $I_1$, $I_2$ are fixed points}
Our remaining task is to check \eqref{eq:assum} for $I_1 \not= I_2$.
In fact, it is not necessary to check \eqref{eq:assum} for all $I_1$,
$I_2$ thanks to the torus action. As above, we take a one-parameter
subgroup $\lambda\colon \C^*\to T$ and consider the limit of
$\lambda(t)^*I_1$, $\lambda(t)^*I_2$ when $t\to 0$. We may assume both
converge to fixed points of the torus action.
If $\lambda(t)^*I_1$ and $\lambda(t)^*I_2$ converge to the same point,
it means that they are contained in a neighborhood of the diagonal for
sufficiently small $t$. This case was treated in the previous
subsection. Thus we may assume that $\lambda(t)^*I_1$ and
$\lambda(t)^*I_2$ converge to different points if $I_1 \neq I_2$.  By
the semicontinuity, we have
\begin{equation*}
\dim \Hom_A(I_1, A/I_2)^G \leq
\dim \Hom_A(\lim_{t\to 0}\lambda(t)^*I_1, 
A/\lim_{t\to 0}\lambda(t)^*I_2))^G.
\end{equation*}
Since the left hand side is bounded by $1$ from below (see
\lemref{lem:H1}), it is enough to show the right hand side is $1$.
Thus we may assume $I_1$ and $I_2$ are different fixed
points of the torus action.

\subsection{Classification of the pair of fixed points}
In the remaining of this paper, we assume $I_1$ and $I_2$ are different fixed
point of the torus action.
By Proposition~\ref{prop:AB}, we have four possibilities:
\begin{align}
& \text{$I_1$ of type (A), $I_2$ of type (A)}, \tag*{(AA)} \\ 
& \text{$I_1$ of type (B), $I_2$ of type (B)}, \tag*{(BB)} \\
& \text{$I_1$ of type (A), $I_2$ of type (B)}, \tag*{(AB)} \\
& \text{$I_1$ of type (B), $I_2$ of type (A)}. \tag*{(BA)}
\end{align}
In order to treat the cases (A) and (B) simultaneously, we write the
exponents of generators of $x^\bullet$, $y^\bullet$, $z^\bullet$ by
$\alpha$, $\beta$, $\gamma$ as before. Namely, $\alpha= a+d-1$
(resp.~$a+d$) in case (A) (resp.~(B)), etc. 
We put ``prime'' on the exponents of the generators for the ideal
$I_2$. Namely, $a'$, $\alpha'$, etc.
In either cases, we have
\begin{equation}\label{eq:relation}
\begin{split}
  a+d-1 \leq \alpha \leq a+d, \quad
  b+e-1 & \leq \beta \leq b+e, \quad
  c+f-1 \leq \gamma \leq c+f, \\
  a'+d'-1 \leq \alpha' \leq a'+d', \quad
  b'+e'-1 & \leq \beta' \leq b'+e', \quad
  c'+f'-1 \leq \gamma' \leq c'+f'
\end{split}
\end{equation}
by \propref{prop:AB}.

\begin{Lemma}\label{lem:yes}
At least one of  three generators $x^\alpha$, $y^\beta$ and $z^\gamma$
of $I_1$ belongs to the ideal $I_2$.
\end{Lemma}

\begin{proof}
We assume $x^\alpha$, $y^\beta$, $z^\gamma \not\in I_2$ and derive a
contradiction.
From this assumption we have
\begin{equation}\label{eq:1}
  \alpha < \alpha',\quad
  \beta < {\beta'},\quad
  \gamma < {\gamma'}.
\end{equation}

Since $x^{a-1}y^{e-1}$ and $z^\gamma$ define the isomorphic
representation of $G$, we have $x^{a-1}y^{e-1} \in I_2$ by
\lemref{lem:rep}, as $z^\gamma\notin I_2$ by the assumption. Then one
of the following holds:
\begin{subequations}
\begin{alignat}{3}
 & a' \leq a-1 &\quad& \text{and} &\quad& e' \leq e-1, \label{eq:1a} \\
 & a' > a-1 &&\text{and}&& {\beta'} \leq e-1, \label{eq:1i} \\
 & \alpha' \leq a-1 &&\text{and} &&e' > e-1. \label{eq:1u} 
\end{alignat}
\end{subequations}

The case \eqref{eq:1i} contradicts with \eqref{eq:1}, \eqref{eq:relation} as 
$\beta < {\beta'} \leq e - 1 < \beta$. Similarly we have a contradiction in the
case \eqref{eq:1u}. Hence we must have \eqref{eq:1a}.
Exchanging $x$, $y$ and $z$ we have 
\begin{equation*}
 a' \leq a-1, \  b' \leq b-1, \ c' \leq c-1, \ 
 d' \leq d-1, \  e' \leq e-1, \ f' \leq f-1. 
\end{equation*}

Combining these with \eqref{eq:relation}, we get 
\begin{equation*}
\alpha' \leq  a' + d' \leq  a+d -2 \leq \alpha-1 < \alpha'.
\end{equation*}
This is a contradiction.
\end{proof}

\begin{Lemma}\label{lem:Hiraku}
Assume
\begin{equation}\label{eq:*}
{\beta'} \leq \beta \quad \text{and} \quad {\gamma'} \leq \gamma.
\end{equation}
If $x^\alpha = x^{\alpha'}$, then $y^bz^f = y^{b'}z^{f'}$. The same holds if
we exchange $x$, $y$ and $z$.
\end{Lemma}

\begin{proof}
We suppose $y^{b-1}z^{f-1}\neq y^{b'-1}z^{f'-1}$ and lead to a
contradiction.
Both $y^{b-1}z^{f-1}$ and $y^{b'-1}z^{f'-1}$ give the isomorphic
representation as $x^\alpha = x^{\alpha'}$. Since 
$y^{b-1}z^{f-1}\notin I_1$, $y^{b'-1}z^{f'-1}\notin I_2$,
\lemref{lem:rep} implies
\begin{equation*}
  y^{b-1}z^{f-1}\in I_2, \quad y^{b'-1}z^{f'-1}\in I_1.
\end{equation*}
From the first condition, we have one of the following:
\begin{subequations}
\begin{alignat}{3}
 & b' \leq b-1 &\quad&\text{and}&\quad & f' \leq f-1, \label{co:H} \\
 & b' > b-1  &&\text{and}&  & \gamma' \leq f-1, \label{co:R}  \\
 & \beta' \leq b-1 &&\text{and}& & f' > f-1. \label{co:K}
\end{alignat}
\end{subequations}
Similarly the second condition implies the one of the following:
\begin{subequations}
\begin{alignat}{3}
 & b \leq b'-1 &\quad&\text{and}&\quad & f \leq f'-1, \label{co:I} \\
 & b > b'-1  &&\text{and}&  & \gamma \leq f'-1, \label{co:A}  \\
 & \beta \leq b'-1 &&\text{and}& & f > f'-1. \label{co:U}
\end{alignat}
\end{subequations}
The condition~\eqref{co:H} contradicts with \eqref{co:I} as
\begin{equation*}
  b \leq b'-1 \leq b-2.
\end{equation*}
The condition~\eqref{co:R} (resp.~\eqref{co:K}) is not compatible with \eqref{co:I} and \eqref{eq:relation} because 
\begin{equation*}
{\gamma'} \leq f-1 \leq f'-2 \leq {\gamma'} -2 \
(\text{resp.}  \ {\beta'} \leq b-1 \leq b'-2 \leq {\beta'} -2).
\end{equation*}
Hence the case \eqref{co:I} does not hold.

And either the condition \eqref{co:A} or \eqref{co:U} contradicts with
the conditions \eqref{eq:relation} and \eqref{eq:*} as
\begin{align}
& \gamma \leq f'-1 \leq {\gamma'} -1 \leq \gamma -1,\tag*{Case \eqref{co:A}:} \\
& \beta \leq b'-1 \leq {\beta'} -1 \leq \beta -1.\tag*{Case \eqref{co:U}:} 
\end{align}
This completes the proof.
\end{proof}

\begin{figure}[htbp]
\begin{center}
\includegraphics{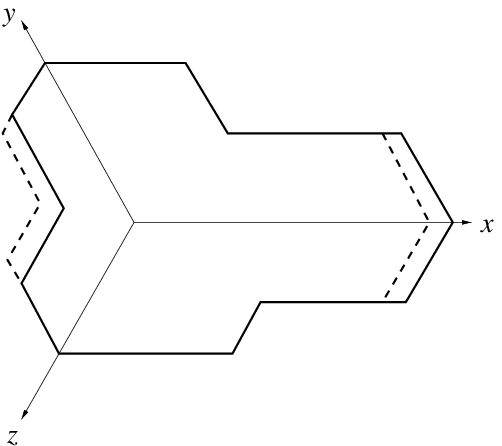}
\caption{Case C}
\label{fig:C}
\end{center}
\end{figure}

\begin{Lemma}\label{lem:no}
\rom{(1)} When the pair $(I_1,I_2)$ is not of the type \rom{(BA)},
at least one of $x^\alpha, y^\beta$ and $z^\gamma$ in $I_1$ does not belong
to $I_2$.

\rom{(2)} When the pair $(I_1,I_2)$ is of the type \rom{(BA)},
one of the following two cases occurs:
\begin{aenume}
\item at least one of $x^\alpha$, $y^\beta$ and $z^\gamma$ in $I_1$ does not
belong to $I_2$,
\item we have
\begin{equation}\label{eq:(4)}
a = a',\ b = b'-1,\ c = c',\
d = d',\ e = e',\   f = f' - 1 
\end{equation}
or the one with $x$, $y$, $z$ exchanged.
(See Figure~\ref{fig:C}.)\end{aenume}
\end{Lemma}

\begin{proof}
Assume that $x^\alpha$, $y^\beta$, $z^\gamma \in I_2$. Then we have
\begin{equation*}
  \alpha' \leq \alpha,\
  {\beta'} \leq \beta,\
  {\gamma'} \leq \gamma.
\end{equation*}

After exchanging $x$, $y$, $z$ if necessary, we may assume either of
the following 3 cases occurs:
\begin{subequations}
\begin{align}
 & \alpha' < \alpha, \ {\beta'} < \beta, \ {\gamma'} \leq\gamma, \label{ca:a}\\
 & \alpha' < \alpha, \ {\beta'} =\beta, \ {\gamma'} =\gamma, \label{ca:u}\\
 & \alpha' = \alpha, \ {\beta'} = \beta, \ {\gamma'} =\gamma \label{ca:e}.
\end{align}
\end{subequations}

In case~\eqref{ca:a}, we have $x^{\alpha'},y^{\beta'}\not\in I_1$, and hence 
$y^{b'-1}z^{f'-1}$, $z^{c'-1}x^{d'-1} \in I_1$ by \lemref{lem:rep}.
From $y^{b'-1}z^{f'-1} \in I_1$, one of the following holds:
\begin{subequations}
\begin{alignat}{3}
 & b \leq b'-1 &\quad&\text{and}&\quad & f \leq f'-1, \label{eq:Yu} \\
 & b > b'-1  &&\text{and}&  & \gamma \leq f'-1, \label{eq:Ka}  \\
 & \beta \leq b'-1 &&\text{and}& & f > f'-1. \label{eq:Ri} 
\end{alignat}
\end{subequations}
The case \eqref{eq:Ka} contradicts with \eqref{ca:a} and
\eqref{eq:relation} for
$$ {\gamma'} \leq \gamma \leq f'-1 < {\gamma'}. $$
Similarly,
\eqref{eq:Ri} contradicts with \eqref{ca:a}. Hence we must have
\eqref{eq:Yu}. From the same consideration for $z^{c'-1}x^{d'-1}$, we
have $c \leq c' - 1$ and $d \leq d' - 1$.
Therefore we obtain 
$$ c+f \leq c' + f' -2 .$$
However this inequality together with \eqref{eq:relation} contradicts
with the assumption ${\gamma'} \leq \gamma$. Thus \eqref{ca:a} is excluded.

Next consider the case~\eqref{ca:u}. By \lemref{lem:Hiraku}, we have
\(
a = a', c = c', d = d', e = e'.
\)
We repeat the same argument as in \eqref{ca:a} for the condition
$x^{\alpha'}\notin I_1$ to get $b \leq b'-1$, $f \leq f' -1$.
Combining with \eqref{eq:relation}, we get 
${\beta'} = \beta \leq b+e \leq b'+e'-1 \leq {\beta'}$.  Hence we must have equalities, i.e.\  $\beta = b+e$, ${\beta'} = b'+ e'-1$ and $b= b'-1$.
In particular, the pair of the ideals $(I_1,I_2)$ is of type (BA).
Similarly we have 
$f=f'-1$ by exchanging 
$y$ and $z$.
Thus we must have \eqref{eq:(4)}.

Finally, consider the case~\eqref{ca:e}. By \lemref{lem:Hiraku}, we
have $a = a'$, $b = b'$, $c = c'$, $d = d'$, $e = e'$, $f = f'$. These
together with \eqref{ca:e} means $I_1 = I_2$
which is excluded from the beginning. Then we proved the lemma.
\end{proof}

\subsection{Division of cases}
Before starting the proof of \thmref{prop:assum} for $I_1\neq I_2$, 
we divide cases by the number of generators $x^\alpha$,
$y^\beta$, $z^\gamma$ belonged in $I_2$. By \lemref{lem:yes}, the number is
either $1$, $2$, or $3$. After exchanging $x$, $y$, $z$,
we have
\begin{align}
& x^\alpha\in I_2,\ y^\beta\not\in I_2,\ z^\gamma\not\in I_2, \tag*{Case A:}\\
& x^\alpha\in I_2,\ y^\beta\in I_2,\ z^\gamma\not\in I_2, \tag*{Case B:}\\
& x^\alpha\in I_2,\ y^\beta \in I_2, \ z^\gamma \in I_2. \tag*{Case C:}
\end{align}

In case A, we have $\alpha' \leq \alpha$, ${\beta'} > \beta$,
${\gamma'} > \gamma$. The last two inequalities implies that
$y^{\beta'}, z^{\gamma'}\in I_1$. Then we do not have $x^{\alpha'}\in
I_1$ since $(I_2, I_1)$ is not of the type
\lemref{lem:no}(2)(b). Hence we must have $x^{\alpha'}\notin I_1$. In
summary, we have
\begin{equation*}
 x^\alpha\in I_2,\ y^\beta\not\in I_2,\ z^\gamma\not\in I_2,\ 
 x^{\alpha'}\notin I_1,\ y^{\beta'}\in I_1,\ z^{\gamma'}\in I_1.
\end{equation*}
Thus we have
\begin{equation}\label{eq:AAA}
\text{Case A} \Longleftrightarrow \alpha > \alpha',\ \beta < {\beta'},\ \gamma < {\gamma'}.
\end{equation}

In case B, we have $\alpha' \leq \alpha$, ${\beta'} \leq \beta$,
${\gamma'} > \gamma$. We have $z^{\gamma'}\in I_1$. We further
separate cases by the number of generators among $x^{\alpha'}$,
$y^{\beta'}$ belonged to $I_1$, i.e.\  $0$, $1$, or $2$. After
exchanging $x$ and $y$ if necessary, we have
\begin{align}
& x^{\alpha'} \in I_1,\ y^{\beta'} \in I_1, \tag*{Case B1:}\\
& x^{\alpha'} \not\in I_1,\ y^{\beta'} \in I_1, \tag*{Case B2:}\\
& x^{\alpha'} \not\in I_1,\ y^{\beta'}\notin I_1. \tag*{Case B3:}
\end{align}
Since we have $\alpha\geq\alpha'$ (resp.\ $\beta\geq {\beta'}$), $x^{\alpha'}\in I_1$
(resp.\ $y^{\beta'}\in I_1$) is possible only when $\alpha = \alpha'$ (resp.\
$\beta={\beta'}$). Thus we have
\begin{align}
\text{Case B1} &\Longleftrightarrow \alpha = \alpha',\ \beta =
{\beta'},\ \gamma <{\gamma'},\\ 
\text{Case B2} &\Longleftrightarrow \alpha > \alpha',\ \beta =
{\beta'},\ \gamma <{\gamma'},\label{eq:BBB}\\ 
\text{Case B3} &\Longleftrightarrow \alpha
> \alpha',\ \beta > {\beta'},\ \gamma <{\gamma'}.
\end{align}

In case C, we have only one possibility by \lemref{lem:no}(2).

In fact, it is not necessary to check $\dim\Hom_A(I_1, A/I_2)^G = 1$
all cases: In the complex~\eqref{eq:koszul}, we have seen that the
second homology $H^2$ is the dual space of the first cohomology $H^1$
of the complex in which $I_1$ and $I_2$ are exchanged (see
\lemref{lem:H1}). If $I_1 \neq I_2$, then $H^0= H^3 = 0$
(\lemref{lem:H0}). Since the Euler characteristic of the complex is
equals to $0$, $H^1 = 0$ implies $H^2 = 0$, and hence $H^1 =0$ of the
complex with $I_1$ and $I_2$ exchanged. Therefore it is sufficient to
show that $H^1=0$ for one of the two pairs $(I_1,I_2)$ and
$(I_2,I_1)$. If we exchange $I_1$ and $I_2$, then the case A and B3,
B1 and C are exchanged. So we only need to consider the cases A, B2,
C.

Let $\Phi\colon I_1\to A/I_2$ be a $G$-equivariant
$A$-homomorphism. First of all, we have $\Phi(xyz) = 0$ as
follows. Since the representation corresponding to $xyz$ is trivial,
we have $\Phi(xyz) = v\bmod I_2$ for some $v\in\C$. Since $I_1\neq
I_2$, there exists $f\in I_1\setminus I_2$. Then
\begin{equation*}
  vf \bmod I_2 = f \Phi(xyz) = \Phi(xyz\, f) = xyz \Phi(f) = 0 \bmod I_2.
\end{equation*}
Since $f\notin I_2$, we must have $v = 0$.

\subsection{Case A} 
First we consider the case A.
\begin{Lemma}\label{clm:A}
{\rm (i)} $x^a y^e\in I_2$ and $z^c x^d \in I_2$.

{\rm (ii)} $y^bz^{\gamma} \not\in I_2$ and $y^{\beta}z^f \not\in I_2$.
\end{Lemma}
\begin{proof}
(i) Since  $z^{\gamma-1} \not\in I_2$ by \eqref{eq:AAA}, we have $x^ay^e \in I_2$
by \lemref{lem:rep}. Exchanging $y$ and $z$, we also have
$z^cx^d\in I_2$.

(ii)
As $x^{\alpha'} \not\in I_1$, we have $y^{b'-1}z^{f'-1}\in I_1$ by 
\lemref{lem:rep}. Then one of the following holds:
\begin{subequations}
\begin{alignat}{3}
 &  b \leq b'-1 &\quad& \text{and} &\quad& f \leq f'-1, \label{eq:A1} \\
 & b > b'-1 &&\text{and}&& \gamma \leq f'-1, \label{eq:A2} \\
 & \beta \leq b'-1 &&\text{and} &&f > f'-1. \label{eq:A3} 
\end{alignat}
\end{subequations}

Now we suppose $y^bz^{\gamma} \in I_2$, then we have one of the following:
\begin{subequations}
\begin{alignat}{3}
 &  b' \leq b &\quad& \text{and} &\quad& f' \leq \gamma , \label{eq:B1} \\
 & b' > b &&\text{and}&& {\gamma'} \leq \gamma, \label{eq:B2} \\
 & {\beta'} \leq b &&\text{and} &&f' > \gamma. \label{eq:B3} 
\end{alignat}
\end{subequations}

The case \eqref{eq:B2} contradicts with \eqref{eq:AAA} because 
$\gamma < {\gamma'} \leq \gamma$.
Similarly we have a contradiction in the case 
\eqref{eq:B3}. Thus we must have \eqref{eq:B1}.
However the case \eqref{eq:B1} is not compatible with any of
\eqref{eq:A1},\eqref{eq:A2},\eqref{eq:A3}:
\begin{align}
& b' \leq b \leq b'-1,\tag*{Case \eqref{eq:A1}:} \\
&f'\leq \gamma \leq f'-1,\tag*{Case \eqref{eq:A2}:} \\
&\beta \leq b'-1 \leq b-1 \leq \beta -1\tag*{Case \eqref{eq:A3}:} 
\end{align}
Therefore we have $y^bz^\gamma \not\in I_2$.
By exchanging $y$ and $z$, we also have 
$y^{\beta} z^f \not\in I_2$. 
\end{proof}

\begin{figure}[htbp]
\begin{center}
\includegraphics{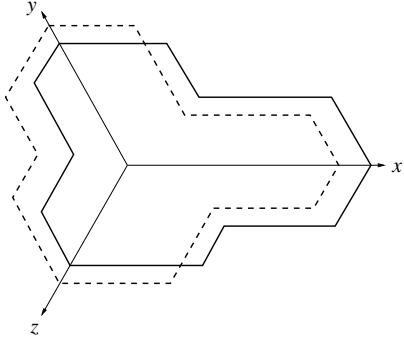}
\caption{Case A: $I_1 =$ solid lines, $I_2 =$ dotted lines}
\label{fig:idealA}
\end{center}
\end{figure}

Thus we can draw junior diagrams for $I_1$ and $I_2$ as
Figure~\ref{fig:idealA}.

Let $\Phi\colon I_1
\rightarrow A/I_2$ be a $G$-equivariant $A$-homomorphism. 
By \lemref{clm:A} and \eqref{eq:AAA}, we can determine the image of
generators of $I_1$ as follows:
\begin{alignat*}{2}
& \Phi (x^\alpha)= p y^{b-1}z^{f-1} \bmod I_2,& \quad 
& \Phi (x^ay^e)=s z^{\gamma-1} \bmod I_2,\\
& \Phi (y^\beta)= q y^\beta \bmod I_2,& \quad 
& \Phi (y^bz^f)=t y^bz^f \bmod  I_2,\\
& \Phi (z^\gamma)= r z^\gamma \bmod I_2,& \quad 
& \Phi (z^cx^d)=u y^{\beta-1} \bmod I_2.
\end{alignat*}
Here we determine the image so that it has the isomorphic
representation of $G$ as the generator.

Since $\Phi (y^\beta z^f) = y^\beta \Phi (z^f) = z^f \Phi (y^\beta)$, we have 
$t y^\beta z^f \bmod I_2 = q y^\beta z^f \bmod I_2$. 
Thus we obtain $t=q$ by \lemref{clm:A}(ii). 
Similarly we have $t=r$ by exchanging $y$ and $z$. 

Consider the image of $z^{\gamma}x^d$. We have
$\Phi (z^{\gamma}x^d)=r z^{\gamma}x^d \bmod I_2=0$ by \linebreak[4]
%
%
\lemref{clm:A}(i).
On the other hand, we have $\Phi (z^{\gamma}x^d)=u z^{\gamma-c}y^{\beta-1} \bmod I_2$
where $z^{\gamma-c}y^{\beta-1} \not\in I_2$ by \lemref{clm:A}(ii).
Therefore we have $u=0$. Exchanging $z$ and $y$, we also have $s=0$.

Next see the image of $x^{\alpha}z^c$. We have $\Phi (x^{\alpha}z^c) =
x^{\alpha-d}\Phi (z^cx^d) \bmod I_2= 0$ because $u=0$.  And we have $\Phi (x^{\alpha}z^c) =
p y^{b-1}z^{c+f-1} \bmod I_2$ where $y^{b-1}z^{c+f-1} \not\in I_2$
by \lemref{clm:A}(ii).
Then we obtain $p=0$.

Therefore we see $\dim \Hom_A(I_1, A/I_2)^G = 1$ in case A.

\subsection{Case B2}
Next we study the case B2.
\begin{Lemma}\label{clm:B1}
{\rm (i)} $y^{\beta} \in I_2$,

{\rm (ii)} $x^a y^e \in I_2$ and $z^c x^d \in I_2$.

{\rm (iii)} $f \leq f'-1$.

\end{Lemma}
\begin{proof}
(i) Obvious from  \eqref{eq:BBB}.

(ii) We have $x^ay^e \in I_2$ by 
$z^{\gamma -1}\not\in I_2$ with \lemref{lem:rep}.
We also have $z^cx^d\in I_2$ since $y^{\beta-1} = y^{\beta'-1}\notin I_2$.

(iii) As $x^{\alpha'}\notin I_1$, we have $y^{b'-1}z^{f'-1}\in I_1$ by
\lemref{lem:rep}. As in the proof of \lemref{clm:A}(ii), we have
\eqref{eq:A1}, \eqref{eq:A2}, or \eqref{eq:A3}. 
But \eqref{eq:A3} contradicts with \eqref{eq:BBB}, as
$\beta \le b'-1 \le \beta'-1 = \beta - 1$.
Thus we have either \eqref{eq:A1} or \eqref{eq:A2}.
In either cases, we have the assertion.
\end{proof}

\begin{Lemma}\label{clm:B2}
{\rm (i)} $y^{\beta -1}z^f \not\in I_2$.

{\rm (ii)} $y^{b-1}z^{\gamma}\not\in I_2$. 
\end{Lemma}
\begin{proof} 
(i) If we assume $y^{\beta -1}z^f \in I_2$, then one of the following holds:
\begin{subequations}
\begin{alignat}{3}
 &  b' \leq \beta -1 &\quad& \text{and} &\quad& f' \leq f, \label{eq:11} \\
 & b' > \beta -1 &&\text{and}&& {\gamma'} \leq f, \label{eq:22} \\
 & {\beta'} \leq \beta -1 &&\text{and} && f' > f. \label{eq:33} 
\end{alignat}
\end{subequations}
\eqref{eq:11} contradicts with \lemref{clm:B1}(iii) as
$f' \leq f \leq f' -1$. \eqref{eq:22} is not compatible with 
\eqref{eq:BBB} because $\beta -1 < b' \leq {\beta'} =\beta$. \eqref{eq:33} 
also contradicts with \eqref{eq:BBB} as $ {\beta'} \leq \beta -1 = {\beta'} -1$.
Thus $y^{\beta -1}z^f \not\in I_2$.

(ii) 
If we assume $y^{b -1}z^{\gamma} \in I_2$, then one of the following holds:
\begin{subequations}
\begin{alignat}{3}
 &  b' \leq b -1 &\quad& \text{and} &\quad& f' \leq \gamma, \label{eq:C1} \\
 & b' > b -1 &&\text{and}&& {\gamma'} \leq \gamma, \label{eq:C2} \\
 & {\beta'} \leq b -1 &&\text{and} && f' > \gamma. \label{eq:C3} 
\end{alignat}
\end{subequations}
The condition \eqref{eq:C1} contradicts with \eqref{eq:A1} because
$b \leq b'-1 \leq b -2$. And \eqref{eq:C1} is not compatible with 
\eqref{eq:A2} as $f' \leq \gamma \leq f'-1$.
Thus the case \eqref{eq:C1} is excluded.
The condition \eqref{eq:C2} (resp.~\eqref{eq:C3}) contradicts with the condition \eqref{eq:BBB} as
$$ {\gamma'} \leq \gamma < {\gamma'} \ (\text{resp.}\  {\beta'} \leq b-1 \leq \beta -1 = {\beta'} -1). $$
Thus the proof completes.
\end{proof}

We study the cases $b\neq \beta$ and $b=\beta $ separately. First assume
$b\neq \beta$.

\begin{Lemma}\label{clm:B3}
 $y^bz^{\gamma}\not\in I_2$. 
\end{Lemma}

\begin{proof} 
Now we assume $y^bz^{\gamma}\in I_2$ and derive a contradiction.
From this assumption we have \eqref{eq:B1}, \eqref{eq:B2} 
or \eqref{eq:B3}. The cases \eqref{eq:B1} and
\eqref{eq:B2} lead to contradictions as in \lemref{clm:A}(ii). 
Thus we have case \eqref{eq:B3}.
However it implies
${\beta'} \leq b \leq \beta ={\beta'}$, hence $\beta = b$.
This contradicts with the assumption.
\end{proof}

\begin{figure}[htbp]
\begin{center}
\includegraphics{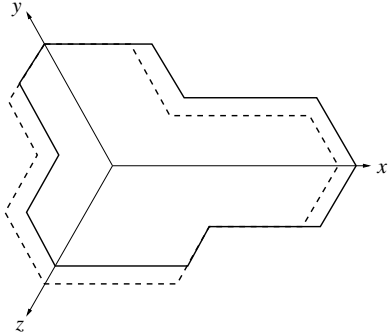}
\caption{Case B2 ($c = c'$, $d = d'$)}
\label{fig:B2a}
\end{center}
\end{figure}

\begin{figure}[htbp]
\begin{center}
\leavevmode
\includegraphics{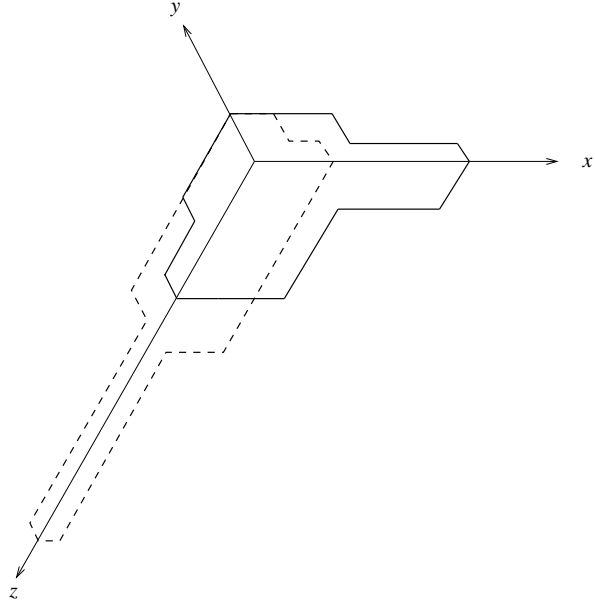}
\caption{B2 ($c\neq c'$ or $d\neq d'$)}
\label{fig:B2b}
\end{center}
\end{figure}

By above, we can draw junior diagrams as Figure~\ref{fig:B2a} 
(if $c = c'$, $d = d'$) or Figure~\ref{fig:B2b} (otherwise).

By above discussions, we can put the images as follows:
\begin{alignat*}{2}
& \Phi (x^\alpha)= p y^{b-1}z^{f-1} \bmod I_2,& \quad 
& \Phi (x^ay^e)=s z^{\gamma-1} \bmod I_2,\\
& \Phi (y^\beta)= q z^{c'-1}x^{d'-1} \bmod I_2,& \quad 
& \Phi (y^bz^f)=t y^bz^f \bmod  I_2,\\
& \Phi (z^\gamma)= r z^\gamma \bmod I_2,& \quad 
& \Phi (z^cx^d)=u y^{\beta-1} \bmod I_2.
\end{alignat*}
Here $\Phi(y^\beta)$ was determined by using $\beta = {\beta'}$.

As in case A, we have $r=t$ by \lemref{clm:B3}. 
We have $u=0$ by Lemmas~\ref{clm:B1}(ii) and \ref{clm:B2}(i) as in case A. We also obtain
$p=0$  by the same discussion in case A with Lemmas~\ref{clm:B1}(ii) and 
\ref{clm:B2}(ii).

Now we consider the image of $y^{\beta}z^f$.  We have
$\Phi(y^{\beta}z^f) = t y^{\beta}z^f \bmod I_2=0$ by
\lemref{clm:B1}(i).  On the other hand, we have $\Phi(y^{\beta}z^f) =
q z^{c'+f-1}x^{d'-1} \bmod I_2$ and $z^{c'+f-1}x^{d'-1}\not\in I_2$
because of \lemref{clm:B1}(iii). Thus we get $q=0$.

For the image of $x^ay^{\beta}$, we have 
$\Phi(x^ay^{\beta}) = x^a \Phi(y^{\beta})=0$ because $q=0$,
while $\Phi(x^ay^{\beta}) = s y^{\beta-e}z^{\gamma-1}\bmod  I_2$ and
$y^{\beta-e}z^{\gamma-1} \not\in I_2$ by \lemref{clm:B3}.
Therefore we obtain $s=0$ and we see the dimension is also one.

Next assume $b =\beta$. Then $I_1$ is of type (A).
In this case, neither \eqref{eq:A1} nor \eqref{eq:A3} holds because
$ \beta = b \leq b' - 1 \leq \beta' -1 = \beta -1 $. Thus the condition
\eqref{eq:A2} always holds, i.e.\ 
\begin{equation}\label{eq:FF}
 b > b' -1 \quad \text{and} \quad \gamma \leq f' -1.
\end{equation} 
The images of the generators are written as follows:
\begin{alignat*}{2}
& \Phi (x^\alpha)= p y^{b-1}z^{f-1} \bmod I_2,& \quad 
& \Phi (x^ay^e)=s z^{\gamma-1} \bmod I_2,\\
& \Phi (y^\beta)= q z^{c'-1}x^{d'-1} \bmod I_2,& \quad 
& \\
& \Phi (z^\gamma)= r z^\gamma \bmod I_2,& \quad 
& \Phi (z^cx^d)=u y^{\beta-1} \bmod I_2.
\end{alignat*}
We have $u=0$ and $p=0$ as in case $b \neq \beta$.

Let us consider the image of $y^\beta z^\gamma$. We have $\Phi
(y^\beta z^\gamma ) = r y^\beta z^\gamma \bmod I_2 =0$. On the other
hand, we have $\Phi (y^\beta z^\gamma ) = q z^{\gamma + c' -1} x^{d' -
1} \bmod I_2$ where $z^{\gamma + c' -1} x^{d' -1} \not\in I_2$ by
\eqref{eq:FF}.  Therefore we have $q=0$. Then we can use same
discussion to obtain $s=0$ as in case $b \neq \beta$.  Thus we have
the assertion in this case.

\subsection{Case C}
In case C, we already see the relation between the exponents of the 
generators of the ideals $I_1$ and $I_2$ in \lemref{lem:no}.

Then we can write as follows:
\begin{alignat*}{2}
& \Phi (x^\alpha)= p y^{b-1}z^{f-1} \bmod I_2,& \quad 
& \Phi (x^ay^e)=s z^{\gamma-1} \bmod I_2,\\
& \Phi (y^\beta)= q z^{c-1}x^{d-1} \bmod I_2,& \quad 
& \Phi (y^bz^f)=t y^bz^f \bmod  I_2,\\
& \Phi (z^\gamma)= r x^{a-1}y^{e-1} \bmod I_2,& \quad 
& \Phi (z^cx^d)=u y^{\beta-1} \bmod I_2.
\end{alignat*}
By similar discussions as in case A and B2, we
obtain $p=q=r=s=u=0$. Then we have the assertion in case C.

Therefore we proved Theorem~\ref{prop:assum}.

\end{document}